\documentclass{amsart}
\usepackage{amssymb}
\usepackage{epsfig}
\usepackage{graphicx}
\numberwithin{equation}{section}
\usepackage[small]{diagrams}

\input xy
\xyoption{all}

\theoremstyle{plain}
\newtheorem{prop}{Proposition}[section]

\newtheorem{theo}[prop]{Theorem}

\newtheorem{lemm}[prop]{Lemma}

\theoremstyle{definition}

\newtheorem{rema}[prop]{Remark}

\newtheorem{exam}[prop]{Example}

\newcommand{\bA}{\mathbb A}
\newcommand{\bC}{\mathbb C}
\newcommand{\bN}{\mathbb N}
\newcommand{\bP}{\mathbb P}
\newcommand{\bQ}{\mathbb Q}
\newcommand{\bZ}{\mathbb Z}
\newcommand{\bR}{\mathbb R}
\newcommand{\bF}{\mathbb F}
\newcommand{\R}{\mathbb R}
\newcommand{\Z}{\mathbb Z}
\newcommand{\bK}{\mathbf K}

\newcommand{\sB}{\mathsf B}
\newcommand{\sH}{\mathsf H}
\newcommand{\sN}{\mathsf N}
\newcommand{\sT}{\mathsf T}
\newcommand{\sZ}{\mathsf Z}

\newcommand{\rk}{\mathrm{rk}}

\newcommand{\ra}{\rightarrow}

\newcommand{\Spec}{{\rm Spec}}
\newcommand{\Pic}{{\rm Pic}}
\newcommand{\Val}{{\rm Val}}
\newcommand{\Hom}{{\rm Hom}}

\newcommand{\mm}{\mathfrak m}

\makeatother
\makeatletter

\author{Sho Tanimoto}
\address{Courant Institute, NYU \\ 251 Mercer str. \\ New York, NY 10012 \\ USA}
\email{tanimoto@cims.nyu.edu}
\author{Yuri Tschinkel}
\address{Courant Institute, NYU \\ 251 Mercer str. \\ New York, NY 10012 \\ USA}
\email{tschinkel@cims.nyu.edu}

\title[Height zeta functions]{Height zeta functions of equivariant compactifications of semi-direct 
products of algebraic groups}

\begin{document}
\date{\today}

\maketitle

\begin{abstract}
We apply the theory of height zeta functions to study the asymptotic distribution of 
rational points of bounded height on projective equivariant compactifications of 
semi-direct products. 
\end{abstract}

\tableofcontents

\section*{Introduction}

Let $X$ be a smooth projective variety over a number field $F$ and $L$ a very ample line bundle on $X$. 
An adelic metrization $\mathcal L = ( L, \| \cdot\|)$ on $L$ induces a height function
$$
\sH_{\mathcal L} \colon X(F)\ra \R_{ >0},
$$
let 
$$
\sN(X^{\circ}, \mathcal L, \sB) :=  \# \{ x\in X^{\circ}(F)\, |\, \sH_{\mathcal L}(x)\le \sB\}, \quad X^{\circ}\subset X, 
$$
be the associated counting function for a subvariety $X^{\circ}$.  
Manin's program, initiated in \cite{fmt} and significantly developed over the last 10 years, 
relates the asymptotic of the counting function
$\sN(X^{\circ},\mathcal L, \sB)$, as $\sB\ra \infty$, for a suitable Zariski open $X^{\circ}\subset X$, 
to global geometric invariants of the underlying variety $X$. By general principles of diophantine geometry, 
such a connection can be expected for varieties with sufficiently positive anticanonical line bundle $-K_X$, 
e.g., for Fano varieties. Manin's conjecture asserts that
\begin{equation}
\label{eqn:main}
\sN(X^{\circ}, -\mathcal K_X, \sB)= c \cdot \sB \log(\sB)^{r-1},
\end{equation}
where $r$ is the rank of the Picard group $\Pic(X)$ of $X$, 
at least over a finite extension of the ground field. 
The constant $c$ admits a conceptual interpretation, its main ingredient
is a Tamagawa-type number introduced by Peyre \cite{peyre95}. 

For recent surveys highlighting different aspects of this program, see, e.g., \cite{t-survey}, 
\cite{chambert-survey}, \cite{browning-survey}, \cite{browning-book}. 

Several approaches to this problem have evolved:
\begin{itemize}
\item passage to (universal) torsors combined with lattice point counts;
\item variants of the circle method;
\item ergodic theory and mixing;
\item height zeta functions and spectral theory on adelic groups.
\end{itemize}

The universal torsor approach has been particularly successful in the treatment of 
del Pezzo surfaces, especially the singular ones. This method works best over 
$\mathbb Q$; applying it to surfaces over more general number fields often presents
insurmountable difficulties, see, e.g., \cite{fouvry}.
Here we will explain the basic principles of the method of height zeta functions
of equivariant compactifications of linear algebraic groups and apply it
to semi-direct products; this method is insensitive to the ground field.  
The spectral expansion of the height zeta function 
involves 1-dimensional as well as infinite-dimensional representations, see Section~\ref{sect:harm}
for details on the spectral theory.  
We show that the main term appearing in the spectral analysis, namely, the term
corresponding to 1-dimensional representations, 
matches precisely the predictions of Manin's conjecture, i.e., has the form \eqref{eqn:main}. 
The analogous result for the universal torsor approach can be found in \cite{peyre-principal} and 
for the circle method applied to universal torsors in \cite{peyre-circle}.

Furthermore, using the tools developed in Section~\ref{sect:harm}, 
we provide new examples of rational surfaces satisfying Manin's conjecture.

\

{\bf Acknowledgments.} 
The second author was partially supported by NSF grants DMS-0739380 and 
0901777.

\section{Geometry}
\label{sect:geom}

In this section, we collect some general 
geometric facts concerning equivariant compactifications of 
solvable linear algebraic groups.  
Here we work over an algebraically closed field of characteristic 0. 

Let $G$ be a linear algebraic group. In dimension 1, the only examples are 
the additive group $\mathbb G_a$ and the multiplicative group $\mathbb G_m$. 
Let
$$
\mathfrak X(G)^{\times}:= \Hom (G, \mathbb G_m),
$$
the group of algebraic characters of $G$. 
In the situations we consider, 
this is a torsion-free $\Z$-module of finite rank
(see \cite[Section 8]{Borel}, for conditions insuring this property).

Let $X$ be a projective 
equivariant compactification of $G$. 
After applying equivariant resolution 
of singularities, if necessary,  
we may assume that $X$ is 
smooth and that the boundary 
$$
X\setminus G= D = \cup_{\iota} D_{\iota},
$$
is a divisor with normal crossings. 
Here $D_{\iota}$ are irreducible components of $D$. 
Let $\Pic^G(X)$ be the group of equivalence classes of $G$-linearized line bundles
on $X$. 

Generally, we will identify divisors, associated line bundles, and 
their classes in $\Pic(X)$, resp. $\Pic^G(X)$.

\begin{prop}
\label{prop:picard}
Let $X$ be a smooth and proper 
equivariant compactification of a connected solvable 
linear algebraic group $G$. Then,
\begin{enumerate}
\item we have an exact sequence
$$
0 \ra \mathfrak X(G)^{\times} \ra \textnormal{Pic}(X)^{G} \ra \textnormal{Pic}(X) \ra 0,
$$
\item $\textnormal{Pic}^G(X) = \oplus_{\iota \in \mathcal I}\bZ D_{\iota}$, and
\item the closed cone of pseudo-effective divisors of $X$ is spanned by the boundary components:
$$
\Lambda_{\textnormal{eff}}(X) = \sum_{\iota \in \mathcal I}\bR_{\geq 0} D_{\iota}.
$$
\end{enumerate}
\end{prop}

\begin{proof}
The first claim follows from the proof of \cite[ Proposition 1.5]{GIT}. 
The crucial point is to show that the 
Picard group of $G$ is trivial. 
As an algebraic variety, a connected solvable group 
is a product of an algebraic torus and an affine space. 
The second assertion holds since 
every finite-dimensional representation of a solvable group has a fixed vector.
For the last statement, see \cite[Theorem 2.5]{HT99}. 
\end{proof}

\begin{prop}
\label{prop:haar}
Let $X$ be a smooth and proper equivariant compactification for 
the left action of a linear algebraic group. Then the right invariant 
top degree differential form $\omega$ satisfies
$$
-\textnormal{div} (\omega) = \sum_{\iota \in \mathcal I} d_{\iota} D_{\iota},
$$
where $d_{\iota} > 0$. The same result holds for the right action and the left invariant form.
\end{prop}
\begin{proof}
See \cite[Theorem 2.7]{HT99}.
\end{proof}

\begin{prop}
\label{prop:contraction}
Let $X$ be a smooth and proper equivariant 
compactification of a connected linear algebraic group. 
Let $f:X \ra Y$ be a birational morphism to a normal 
projective variety $Y$. Then $Y$ is an equivariant compactification of 
$G$ such that the contraction map $f$ is a $G$-morphism.
\end{prop}
\begin{proof}
Choose an embedding 
$Y \hookrightarrow \bP^N$, and let $L$ be the pull back of 
$\mathcal O(1)$ on $X$. Since $Y$ is normal, Zariski's main theorem 
implies that the image of the complete linear series 
$|L|$ is isomorphic to $Y$. 
This linear series carries a $G$-linearization \cite[Corollary 1.6]{GIT}. 
Now we apply the same argument as in the proof of 
\cite[Corollary 2.4]{HT99} to this 
linear series, and our assertion follows. 
\end{proof}

The simplest solvable groups are $\mathbb G_a$ and $\mathbb G_m$, as well as their products.  
New examples arise as semi-direct products. For example, let 
$$
\begin{array}{rcl}
\varphi_d : \mathbb G_m &  \ra  & \mathbb G_m =  \rm{GL}_1, \\
    a & \mapsto & a^d
\end{array}
$$ 
and put
$$
G_d := \mathbb G_a \rtimes_{\varphi_d} \mathbb G_m,
$$
where the group law is given by
$$
(x, a) \cdot (y, b) = (x + \varphi_d(a)y, ab).
$$
It is easy to see that $G_d\simeq G_{-d}$.

One of the central themes in birational geometry is the problem of 
classification of algebraic varieties. 
The classification of $G$-varieties, i.e., varieties with $G$-actions, 
is already a formidable task. 
The theory of toric varieties, i.e., equivariant compactifications of $G=\mathbb G_m^n$, 
is very rich, and provides a testing ground for many conjectures in 
algebraic and arithmetic geometry. 
See \cite{HT99} for first steps towards a classification of equivariant compactifications
of $G=\mathbb G_a^n$, as well as \cite{sharo}, \cite{ar}, \cite{ar2} 
for further results in this direction. 

Much less is known concerning equivariant compactifications of other solvable groups; indeed,  
classifying equivariant compactifications of $G_d$ is already an  
interesting open question. We now collect several results illustrating 
specific phenomena connected with noncommutativity of $G_d$ and with the necessity 
to distinguish actions on the left, on the right, or on both sides.
These play a  role in the analysis of height zeta functions in following sections.
First of all, we have

\begin{lemm}
\label{lemm:bi}
Let $X$ be a biequivariant compactification of a semi-direct product 
$G \rtimes H$ of linear algebraic groups. 
Then $X$ is a one-sided (left- or right-)
equivariant compactification of $G\times H$. 
\end{lemm}

\begin{proof}
Fix one section $s: H \ra G \rtimes H$. Define a left action by
$$
(g, h)\cdot x = g\cdot x \cdot s(h)^{-1},
$$
for any $g \in G$, $h \in H$, and $x \in X$.
\end{proof}

In particular, there is no need to invoke noncommutative harmonic analysis in the treatment of
height zeta functions of biequivariant compactifications 
of general solvable groups since such groups are semi-direct products of 
tori with unipotent groups and the 
lemma reduces the problem to a one-sided action of the {\em direct} product.
Height zeta functions of direct products of additive groups and tori
can be treated by combining the methods of \cite{BT98} and \cite{BT-general} with \cite{VecIII}, 
see Theorem~\ref{theo:mainterm}. However, Manin's conjectures are still open 
for one-sided actions of unipotent groups, even for the Heisenberg group.

The next observation is that the projective plane $\bP^2$ is an equivariant 
compactification of $G_d$, for any $d$. 
Indeed, the embedding 
$$
(x,a) \mapsto (a:x:1) \in \bP^2
$$
defines a left-sided equivariant compactification, with boundary a union of two lines.  
In contrast, we have
 
\begin{prop}
\label{prop:plane_semidirectproduct}
If $d \neq 1, 0,$ or  $-1$, then $\bP^2$ is not a biequivariant compactification of $G_d$.
\end{prop}

\begin{proof}
Assume otherwise. Let $D_1$ and $D_2$ be the two irreducible boundary components. 
Since $\mathcal O(K_{\bP^2}) \cong \mathcal O(-3)$, either both components $D_1$ and $D_2$ 
are lines or one of them is a line and the other a conic. 
Let $\omega$ be a right invariant top degree differential form. Then 
$\omega/\varphi_d(a)$ is a left invariant differential form. 
If one of $D_1$ and $D_2$ is a conic, 
then the divisor of $\omega$ takes the form
$$
-\text{div}(\omega) = - \text{div}(\omega/\varphi_d(a)) = D_1 + D_2,
$$ 
but this is a contradiction. 
If $D_1$ and $D_2$ are lines, then without loss of generality, we can assume that 
$$
-\text{div}(\omega) = 2D_1 + D_2 \quad \text{ and } \quad 
- \text{div}(\omega/\varphi_d(a)) = D_1 + 2D_2.
$$ 
However, $\text{div}(a)$ is a multiple of $D_1 - D_2$, which is also a contradiction.
\end{proof}

Combining this result with Proposition~\ref{prop:contraction}, 
we conclude that a del Pezzo surface 
is not a biequivariant compactification of 
$G_d$, for $d \neq 1, 0, \text{or}, -1$.
Another sample result in this direction is:

\begin{prop}
\label{prop:A3A1}
Let $S$ be the singular quartic del Pezzo surface of type $A_3 + A_1$ defined by
$$
x_0^2 + x_0x_3 +x_2x_4 = x_1x_3 - x_2^2 = 0
$$
Then $S$ is a one-sided equivariant compactification of $G_1$, 
but not a biequivariant compactification of $G_d$ if $d \neq 0$.
\end{prop}

\begin{proof}
For the first assertion, see \cite[Section 5]{DL10}.
Assume that $S$ is a biequivariant compactification of $G_d$. 
Let $\pi :\widetilde{S}\ra S$ be its minimal desingularization. 
Then $\widetilde{S}$ is also a biequivariant compactification of $G_d$. 
It has three $(-1)$-curves $L_1$, $L_2$, and $L_3$, which are the 
strict transforms of 
$$
\{x_0=x_1=x_2=0\}, \quad \{x_0+x_3=x_1=x_2=0\}, \quad \text{ and } \quad \{x_0=x_2=x_3=0\},
$$ 
respectively, and has four $(-2)$-curves $R_1$, $R_2$, $R_3$, and $R_4$. 
The nonzero intersection numbers are given by:
$$
L_1.R_1 = L_2.R_1 = R_1.R_2= R_2.R_3=R_3.L_3 =L_3.R_4 =1.
$$
Since the cone of curves is generated by the components of the boundary, 
these negative curves must be in the boundary because each generates an extremal ray. Since the 
Picard group of $\widetilde{S}$ has rank six, the number of boundary components
is seven. Thus, the boundary is equal to the union of these negative curves.

Let $f :\widetilde{S} \ra \bP^2$ be the birational morphism which contracts 
$L_1$, $L_2$, $L_3$, $R_2$, and $R_3$. 
This induces a biequivariant compactification on $\bP^2$. 
The birational map $f\circ \pi^{-1}: S \dashrightarrow \bP^2$ is given by 
$S \ni (x_0:x_1:x_2:x_3:x_4) \mapsto (x_2:x_0:x_3) \in \bP^2$. 
The images of $R_1$ and $R_4$ are $\{y_0=0\}$ and $\{y_2=0\}$ and 
we denote them by $D_0$ and $D_2$, respectively. 
The images of $L_1$ and $L_2$ are $(0:0:1)$ and $(0:1:-1)$, respectively; 
so that the induced group action on $\bP^2$ must fix 
$(0:0:1)$, $(0:1:-1)$, and $D_0\cap D_2 = (0:1:0)$. 
Thus, the group action must fix the line 
$D_0$, and this fact implies that all left and right invariant vector 
fields vanish along $D_0$. It follows that 
$$
-\text{div}(\omega) = - \text{div}(\omega/\varphi_d(a)) = 2D_0 + D_2,
$$ 
which contradicts $d \neq 0$.
\end{proof}


\begin{exam}
\label{exam:Hirz}
Let $l \geq d \geq 0$. The Hirzebruch surface
$\mathbb F_l = \bP_{\bP^1}((\mathcal O\oplus \mathcal O(l))^*)$ 
is a biequivariant compactification of $G_d$.
Indeed, we may take the embedding
$$
\begin{array}{rcc}
G_d  & \hookrightarrow & \mathbb F_l\\
(x,a)&  \mapsto        & ((a:1), [1 \oplus x\sigma_1^l]),
\end{array}
$$
where $\sigma_1$ is a section of the line bundle $\mathcal O(1)$ on $\mathbb P^1$ such that
$$
\text{div}(\sigma_1) = (1:0).
$$
Let $\pi : \bF_l \rightarrow \bP^1$ be the $\bP^1$-fibration. The right action is given by
$$
((x_0:x_1), [y_0 \oplus y_1\sigma_1^l]) \mapsto ((ax_0:x_1), [y_0 \oplus 
(y_1 + (x_0/x_1)^dxy_0)\sigma_1^l]),
$$
on $\pi^{-1}(U_0 = \bP^1 \setminus\{(1:0) \})$ 
and
$$
((x_0:x_1), [y_0 \oplus y_1\sigma_0^l]) 
\mapsto 
((ax_0:x_1), [a^ly_0 \oplus (y_1 + (x_1 /x_0)^{l-d}xy_0)\sigma_0^l]),
$$
on $U_1 = \pi^{-1}(\bP^1 \setminus\{(0:1) \})$. 
Similarly, one defines the left action.
The boundary consists of three components: 
two fibers $f_0 = \pi^{-1}((0:1))$, $f_1 = \pi^{-1}((0:1))$ 
and the special section $D$ characterized by 
$D^2 = -l$.

\end{exam}

\begin{exam}
\label{exam:newrationals}
Consider the right actions in 
Examples \ref{exam:Hirz}. When $l > d >0$, 
these actions fix the fiber $f_0$ and act multiplicatively, i.e., with two fixed points,  
on the fiber $f_1$. Let $X$ be the blow up of two points (or more) on $f_0$ 
and of one fixed point $P$ on $f_1\setminus D$. 
Then $X$ is an equivariant compactification of $G_d$ which is 
neither a toric variety nor a $\mathbb G_a^2$-variety. 
Indeed, there are no equivariant compactifications of 
$\mathbb G_m^2$ on $\bF_l$ fixing $f_0$, so $X$ cannot be toric. 
Also, if $X$ were a $\mathbb G_a^2$-variety, we would obtain an 
induced $\mathbb G_a^2$-action on $\bF_l$ fixing $f_0$ and $P$. 
However,  the boundary consists of two irreducible 
components and must contain $f_0$, $D$, and $P$ 
because $D$ is a negative curve. This is a contradiction.

For $l=2$ and $d=1$, blowing up two points on $f_0$ 
we obtain a quintic del Pezzo surface with an $A_2$ singularity. 
Manin's conjecture for this surface is proved in \cite{derenthal-quintic}. 

In Section~\ref{sect:general}, we prove  Manin's conjecture for $X$ 
with $l \geq 3$.
\end{exam}

\section{Height zeta functions}
\label{sect:height-zeta}

Let $F$ be a number field, $\mathfrak o_F$ its ring of integers, 
and $\Val_F$ the set of equivalence classes of valuations of $F$. 
For $v\in \Val_F$ let $F_v$ be the 
completions of $F$ with respect to $v$, 
for nonarchimedean $v$, let 
$\mathfrak o_v$ be the corresponding ring of integers and 
$\mm_v$ the maximal ideal. Let $\bA=\bA_F$ be the adele ring of $F$. 

Let $X$ be a smooth and projective right-sided equivariant 
compactification of a split connected solvable 
linear algebraic group $G$ over $F$, i.e., 
the toric part $T$ of $G$ is isomorphic to $\mathbb G_m^n$.   
Moreover, we assume that the boundary 
$D = \cup_{\iota \in \mathcal I} D_\iota$ 
consists of geometrically irreducible components meeting transversely. 
We are interested in the asymptotic distribution of rational points of bounded height 
on $X^{\circ}=G\subset X$, with respect to adelically metrized ample line bundles 
$\mathcal L = (L, (\|\cdot \|_{\bA}))$  on $X$. 
We now recall the method of height zeta functions; see \cite[Section 6]{t-survey} for 
more details and examples.

\

\noindent
{\em Step 1.} Define an adelic height pairing
$$
\sH \colon   \Pic^G(X)_{\bC} \times  G(\bA_F) \ra \bC,
$$
whose restriction to 
$$
\sH \colon \Pic^G(X) \times  G(F) \ra \R_{\ge 0},
$$
descends to a {\em height system} on $\Pic(X)$ (see \cite[Definition 2.5.2]{peyre-principal}). 
This means that the restriction of $\sH$ to an $L\in \Pic^G(X)$
defines a Weil height corresponding to some adelic metrization of $L\in \Pic^G(X)$, 
and that it does not depend on the choice of a $G$-linearization on $L$. 
Such a pairing appeared in \cite{BT-aniso} in the context of toric varieties, the extension 
to general solvable groups is straightforward. 

Concretely, by Proposition~\ref{prop:picard}, we know that 
$\Pic^G(X)$ is generated by boundary components $D_{\iota}$, for $\iota \in \mathcal I$. 
The $v$-adic analytic manifold 
$X(F_v)$ admits a ``partition of unity'', i.e., a decomposition into charts
$X_{I,v}$, labeled by $I\subseteq \mathcal I$, such that in each chart
the local height function takes the form
$$
\sH_v({\mathbf s}, x_v)= \phi(x_v)\cdot \prod_{\iota\in I} |x_{\iota,v}|_v^{s_{\iota}} , 
$$
where for each $\iota \in I$,  $x_{\iota}$ is the local coordinate of $D_{\iota}$ 
in this chart,
$$
\mathbf s = \sum_{\iota \in \mathcal I} s_{\iota}D_{\iota},
$$
and $\phi$ is a bounded function, equal to 1 for almost all $v$.  
Note that, locally, the height function 
$$
\sH_{\iota,v}(x_v):= |x_{\iota,v}|_v
$$
is simply the $v$-adic distance to the boundary component $D_{\iota}$. 
To visualize $X_{I,v}$ (for almost all $v$) 
consider the partition induced by 
$$
X(F_v) = X(\mathfrak o_v) \stackrel{\rho}{\longrightarrow} 
\sqcup_{I\subset \mathcal I} X_I^{\circ} (\mathbb F_q), 
$$
where
$$
X_I := \cup_{\iota \in I} D_I, \quad  X_I^{\circ}:=X_I\setminus \cup_{I'\supsetneq I} X_{I'},
$$  
is the stratification of the boundary and $\rho$ is the reduction map; by 
convention $X_{\emptyset}=G$.  
Then $X_{I,v}$ is the preimage of $X_I^{\circ}(\mathbb F_q)$ in $X(F_v)$, and in particular, 
$X_{\emptyset, v} = G(\mathfrak o_v)$, for almost all $v$.  

Since the action of $G$ lifts to integral models of $G$, $X$, and $L$, 
the nonarchimedean local height pairings are invariant with respect to a compact subgroup 
$\bK_v \subset G(F_v)$, which is $G(\mathfrak o_v)$, for almost all $v$.

\

{\em Step 2.}
The height zeta function 
$$
\sZ(\mathbf s, g) := \sum_{\gamma\in G(F)} \sH(\mathbf s, \gamma g)^{-1},
$$
converges absolutely to a holomorphic function, for $\Re(\mathbf s)$ sufficiently large,
and defines a continuous function in 
$\mathsf L^1(G(F)\backslash G(\bA_F)) \cap \mathsf L^2(G(F)\backslash G(\bA_F))$.
Formally, we have the spectral expansion
\begin{equation}
\label{eqn:spectral}
\sZ(\mathbf s, g) = \sum_{\pi} \sZ_{\pi}(\mathbf s, g), 
\end{equation}
where the ``sum'' is over irreducible unitary representations occurring in 
the right regular representation of $G(\bA_F)$ in  
$\mathsf L^2(G(F)\backslash G(\bA_F))$. The invariance of the global height pairing 
under the action of a compact subgroup $\bK\subset G(\bA_F)$, on the side of the action, 
insures that $\sZ_{\pi}$ are in $\mathsf L^2(G(F)\backslash G(\bA_F))^{\bK}$.

\

{\em Step 3.}
Ideally, we would like to obtain a meromorphic continuation of $\sZ$ to a 
tube domain 
$$
\sT_{\Omega} = \Omega  + i\, \Pic(X)_{\R}\subset \Pic(X)_{\bC},  
$$
where $\Omega\subset \Pic(X)_{\R}$ 
is an open neighborhood of the 
anticanonical class $-K_X$. 
It is expected that $\sZ$ is holomorphic for 
$$
\Re(\mathbf s)\in -K_X+ \Lambda_{\rm eff}^{\circ}(X)
$$
and that the polar set of the shifted height zeta function
$\sZ(\mathbf s- K_X, g)$ is the same as that of
\begin{equation}
\label{eqn:laplace}
\mathcal X_{\Lambda_{\rm eff}(X)}(\mathbf s) := 
\int_{\Lambda^*_{\rm eff}(X)} e^{-\langle \mathbf s, \mathbf y\rangle} \mathrm d y,
\end{equation}
the Laplace transform of the set-theoretic characteristic function of
the dual cone $\Lambda_{\rm eff}(X)^*\subset \Pic(X)_{\R}^*$. Here the Lebesgue measure
$\mathrm d y$ is normalized by the dual lattice $\Pic(X)^*\subset \Pic(X)^*_{\R}$. 
In particular, for 
$$
\kappa = -K_X= \sum_{\iota} \kappa_{\iota} D_{\iota},
$$
the restriction of the height zeta function $\sZ(\mathbf s, {\rm id})$ to the one-parameter
zeta function $\sZ(s\kappa,{\rm id})$ should be holomorphic for $\Re(s) >1$, 
admit a meromorphic continuation to $\Re(s)>1-\epsilon$, for some $\epsilon >0$, 
with a unique pole at $s=1$, of order $r=\rk\,\Pic(X)$. 
Furthermore, is is desirable to have some growth estimates in vertical strips. 
In this case, a Tauberian theorem implies
Manin's conjecture \eqref{eqn:main} for the counting function; 
the quality of the error term depends on the growth rate in vertical strips. 
Finally, the leading constant at the pole of $\sZ(s\kappa,{\rm id})$ is essentially
the Tamagawa-type number defined by Peyre. 
We will refer to this by saying that
the height zeta function $\sZ$ satisfies Manin's conjecture;
a precise definition of this class of functions can be found in 
\cite[Section 3.1]{CT-twisted}.

This strategy has worked well and lead to a proof of Manin's conjecture
for the following varieties:
\begin{itemize}
\item toric varieties \cite{BT-aniso}, \cite{BT98}, \cite{BT-general};
\item equivariant compactifications of additive groups $\mathbb{G}_a^n$ \cite{VecIII}; 
\item equivariant compactifications of unipotent groups \cite{shalika-t}, \cite{shalika-t-general};
\item wonderful compactifications of semi-simple groups of adjoint type \cite{shalika-tt}. 
\end{itemize}
Moreover, applications of Langlands' theory of Eisenstein series
allowed to prove Manin's conjecture for flag varieties  \cite{fmt}, 
their twisted products \cite{strauch}, and 
horospherical varieties \cite{strauch-t}, \cite{CT-twisted}.

The analysis of the spectral expansion \eqref{eqn:spectral} is easier when 
every automorphic representation $\pi$ is 1-dimensional, i.e., when $G$ is 
abelian: $G=\mathbb G_a^n$ or $G=T$, an algebraic torus.  
In these cases, \eqref{eqn:spectral} is simply 
the Fourier expansion of the height zeta function and we have, at least formally, 
\begin{equation}
\label{eqn:abelian}
\sZ(\mathbf s, {\rm id})= \int \widehat{\sH}(\mathbf s, \chi) \mathrm d \chi,
\end{equation} 
where
\begin{equation}
\label{eqn:fourier-tr}
\hat{\sH}(\mathbf s, \chi) = \int_{G(\bA_F)} \sH(\mathbf s, g)^{-1} \bar{\chi}(g) \mathrm dg, 
\end{equation}
is the Fourier transform of the height function, 
$\chi$ is a character of $G(F)\backslash G(\bA_F)$, and 
$\mathrm d\chi$ an appropriate measure on the space of automorphic characters. 
For $G=\mathbb G_a^n$, the space of automorphic characters is $G(F)$ itself,
for $G$ an algebraic torus it is (noncanonically) 
$\mathfrak X(G)^{\times}_{\R}\times \mathcal U_G$, 
where $\mathcal U_G$ is a discrete group.

The $v$-adic integration technique developed by Denef and Loeser  
(see \cite{denef-loeser}, \cite{denef-loeser2}, and \cite{denef-loeser3})
allows to compute local Fourier transforms of height functions, in particular, 
for the trivial character $\chi=1$ and almost all $v$ we obtain
$$
\widehat{\sH}_v(\mathbf s, 1) = \int_{G(F_v)} \sH(\mathbf s, g)^{-1} \mathrm dg = 
\tau_v(G)^{-1}
\left(\sum_{I\subset \mathcal I}
\frac{\#X_I^{\circ}(\mathbb F_q)}{q^{{\rm dim}(X)}}
\prod_{\iota\in I} \frac{q-1}{q^{s_{\iota}-\kappa_{\iota} +1}-1}\right),
$$
where $X_I$ are strata of the stratification described in Step 1 and 
$\tau_v(G)$ is the local Tamagawa number of $G$, 
$$
\tau_v(G)= \frac{\# G(\mathbb F_q)}{q^{{\rm dim}(G)}}.
$$
Such height integrals are geometric versions of 
Igusa's integrals; a comprehensive theory in the analytic and adelic setting 
can be found in \cite{volume}.

The computation of Fourier transforms at nontrivial characters
requires a finer partition of $X(F_v)$ which takes into account possible zeroes
of the {\em phase} of the character in $G(F_v)$; see \cite[Section 10]{VecIII} for the
the additive case and \cite[Section 2]{BT-aniso} for the toric case. 
The result is that in the neighborhood of 
$$
\kappa = \sum_{\iota} \kappa_{\iota} D_{\iota}\in \Pic^G(X),
$$ 
the Fourier transform is regularized as follows
$$
\widehat{\sH}(\mathbf s, \chi) = \begin{cases} 
\prod_{v\notin S(\chi)} 
\prod_{\iota \in \mathcal I(\chi)} \zeta_{F,v}(s_{\iota}-\kappa_{\iota}+1) 
\prod_{v\in S(\chi)} \phi_v(\mathbf s, \chi) & G=\mathbb G_a^n, \\
\prod_{v\notin S(\chi)} 
\prod_{\iota \in \mathcal I} \mathsf L_{F,v}(s_{\iota}-\kappa_{\iota}+1 +im(\chi),\chi_u) 
\prod_{v\in S(\chi)} \phi_v(\mathbf s, \chi) & G=T ,
\end{cases}
$$ 
where 
\begin{itemize}
\item  $\mathcal I(\chi)\subsetneq \mathcal I$;  
\item $S(\chi)$ is a finite set of places, which, in general, depends on $\chi$;
\item  $\zeta_{F,v}$ is a local factor of the Dedekind zeta function of $F$ and 
$\mathsf L_{F,v}$ a local factor of a Hecke $\mathsf L$-function;  
\item $m(\chi)$ is the ``coordinate'' of the automorphic character $\chi$ 
of $G=T$ under the embedding $\mathfrak X(G)^{\times}_{\R}\hookrightarrow \Pic^G(X)_{\R}$
in the exact sequence (1) in Proposition~\ref{prop:picard}
and $\chi_u$ is the ``discrete'' component of $\chi$;
\item  
and  $\phi_v(\mathbf s, \chi)$ is a function which is holomorphic and bounded.  
\end{itemize}
In particular, each 
$\widehat{\sH}(\mathbf s, \chi)$ admits a meromorphic continuation as desired and we
can control the poles of each term.  Moreover, at archimedean places we may use
integration by parts with respect to vector fields in the universal enveloping algebra
of the corresponding real of complex group to derive bounds in terms of the ``phase'' 
of the occurring oscillatory integrals, i.e., in terms of ``coordinates'' of $\chi$.

So far, we have not used the fact that $X$ is an equivariant compactification of $G$. 
Only at this stage do we see that the
$\bK$-invariance of the height is an important, in fact, crucial, property that
allows to establish uniform convergence of the right side of the expansion 
\eqref{eqn:spectral}; it insures that 
$$
\widehat{\sH}(\mathbf s, \chi)=0,
$$
for all $\chi$ which are {\em nontrivial} on $\bK$. 
For $G=\mathbb G_a^n$ this means that the trivial representation
is {\em isolated} and that the integral on the right side of Equation~\eqref{eqn:abelian}
is in fact a {\em sum} over a {\em lattice} of integral points in $G(F)$.
Note that Manin's conjecture {\em fails} for nonequivariant compactifications of the affine space,
there are counterexamples already in dimension three \cite{bt-exam}. 
The analytic method described above fails precisely because we cannot 
insure the convergence on the Fourier expansion.  

A similar effect occurs in the noncommutative setting;
one-sided actions {\em do not} guarantee bi-$\bK$-invariance of the height, 
in contrast with the abelian case. 
Analytically, this translates into subtle convergence issues 
of the spectral expansion, in particular, for infinite-dimensional representation.

\begin{theo}
\label{theo:mainterm}
Let $G$ be an extension of an algebraic torus $T$ by a unipotent group $N$ such that 
$[G,G]=N$ over a number field $F$. 
Let $X$ be an equivariant compactification of $G$ over $F$
and 
$$
\sZ(\mathbf s, g) = \sum_{\gamma\in G(F)} \sH(\mathbf s, \gamma g)^{-1},
$$
the height zeta function with respect to an adelic height pairing as in Step 1. 
Let 
$$
\sZ_0(\mathbf s, g) = \int \sZ_{\chi} (\mathbf s, g) \, \mathrm d\chi,
$$
be the integral over all 1-dimensional automorphic representations of $G(\bA_F)$
occurring in the spectral expansion \eqref{eqn:spectral}. Then 
$\sZ_0$ satisfies Manin's conjecture. 
\end{theo}

\begin{proof}
Let 
$$
1\ra N\ra G\ra T\ra 1,
$$
the defining extension. 
One-dimensional automorphic representations of $G(\bA_F)$ are precisely those
which are trivial on $N(\bA_F)$, i.e., these are automorphic characters of $T$. 
The $\bK$-invariance of the height (on one side) insures that only 
unramified characters, i.e., $\bK_T$-invariant characters contribute to
the spectral expansion of $\sZ_0$.  

Let $M = \mathfrak X(G)^\times$ be the group of 
algebraic characters. 
We have
\begin{align*}
\sZ_0(\mathbf s, \text{id}) &= \int_{M_\bR\times \mathcal U_T} 
\int_{G(F)\backslash G(\bA_F)}  \sZ (\mathbf s, g) \bar{\chi}(g) \, \mathrm d g \mathrm d \chi\\
&= \int_{M_\bR\times \mathcal U_T} \int_{G(\bA_F)}  \sH (\mathbf s, g)^{-1} \bar{\chi}(g) \,\mathrm d g \mathrm d \chi\\
& = \int_{M_\bR} \mathsf  F(\mathbf s+im(\chi))\, \mathrm dm,  
\end{align*}
where 
$$
\mathsf F(\mathbf s):= \sum_{\chi\in \mathcal U_T}  \widehat{\sH}(\mathbf s, \chi_u). 
$$
Computations of local Fourier transforms explained above show 
that $\mathsf F$ can be regularized as follows:
$$
\mathsf F(\mathbf s)=\prod_{\iota\in \mathcal I} \zeta_{F}(s_{\iota}-\kappa_{\iota}+1) \cdot 
\mathsf F_{\infty}(\mathbf s),
$$
where $\mathsf F_{\infty}$ is holomorphic for $\Re(s_{\iota}) -\kappa_{\iota} >-\epsilon$, 
for some $\epsilon >0$,  with growth control in vertical strips. 
Now we have placed ourselves into the 
situation considered in \cite[Section 3]{CT-twisted}:
Theorem 3.1.14 establishes analytic properties of integrals
$$
\int_{M_\bR} 
\frac{1}{\prod_{\iota\in \mathcal I} (s_{\iota} - \kappa_{\iota} + im_{\iota})} \cdot 
\mathsf F_{\infty}(\mathbf s +i m) \,
\mathrm dm,
$$
where the image of 
$\iota \colon M_{\R}\hookrightarrow \bR^{\#\mathcal I}$ 
intersects the simplicial cone $\bR^{\#\mathcal I}_{\ge 0}$
only in the origin. 
The main result is that the analytic properties of such integrals match those of 
the $\mathcal X$-function \eqref{eqn:laplace}
of the {\em image} cone under the projection

\

\centerline{
\xymatrix{
0\ar[r] &  M_{\R}\ar@{=}[d]  \ar[r]^{\iota} &  \bR^{\#\mathcal I} \ar@{=}[d]\ar[r]^{\!\!\!\!\pi} &  
 \bR^{\#\mathcal I - \dim(M)} \ar@{=}[d]\ar[r] & 0 \\
0 \ar[r] & \mathfrak X(G)^{\times}_{\R}   \ar[r]  &  \Pic^G(X)_{\R} \ar[r] & \Pic(X)_{\R} \ar[r] & 0; 
}
}

\

\noindent
according to Proposition~\ref{prop:picard}, the image of the simplicial cone 
$\bR^{\#\mathcal I}_{\ge 0}$ under $\pi$ is precisely $\Lambda_{\rm eff}(X)\subset \Pic(X)_{\R}$. 
\end{proof}

\section{Harmonic analysis}
\label{sect:harm}

In this section we study the local and adelic representation theory of
$$
G := \mathbb G_a \rtimes_{\varphi} \mathbb G_m,
$$
an extension of $T:=\mathbb G_m$ by $N:=\mathbb G_a$ via a homomorphism
$\varphi : \mathbb G_m \ra \rm{GL}_1$. 
The group law given by
$$
(x, a) \cdot (y, b) = (x + \varphi(a)y, ab).
$$
We fix the standard Haar measures 
$$
\mathrm dx=\prod_v \mathrm dx_v\quad  \text{and} \quad \mathrm da^\times = \prod_v \mathrm da_v^\times,
$$
on $N(\bA_F)$ and $T(\bA_F)$.  
Note that $G(\bA_F)$ is not unimodular; 
$\mathrm d g:=\mathrm dx \mathrm da^\times$ is a right invariant measure on $G(\bA_F)$
and $\mathrm d g/\varphi(a)$ is a left invariant measure.

Let $\varrho$ be the right regular unitary representation of $G(\bA_F)$ on the Hilbert space:
$$
\mathcal{H}:=\mathsf L^2(G(F)\backslash G(\bA_F), \mathrm d g).
$$
We now discuss the decomposition of $\mathcal{H}$ into irreducible representations. 
Let 
$$
\psi=\prod_v \psi_v : \bA_F \ra \mathbb S^1,
$$ 
be the standard automorphic character and 
$\psi_n$  the character defined by 
$$
x \ra \psi(n x),
$$ 
for $n \in F^{\times}$. 
Let
$$
W:=\ker (\varphi : F^{\times} \ra F^{\times}),
$$ 
and
$$
\pi_{n}:=\text{Ind}_{N(\bA_F) \times W}^{G(\bA_F)}(\psi_n),
$$
for $n \in F^{\times}$.
The following proposition we learned from J. Shalika \cite{shalika-notes}.

\begin{prop}
\label{prop:parameter}
Irreducible automorphic representations, i.e., irreducible unitary 
representations occurring in $\mathcal H = \mathsf L^2(G(F)\backslash G(\mathbb A_F))$,
are parametrized as follows:
$$
\mathcal H = 
\mathsf L^2(T(F)\backslash T(\bA_F)) \oplus 
\widehat{\bigoplus}_{n\in \left(F^{\times}/\varphi(F^{\times})\right)} \pi_{n},
$$
\end{prop}

\begin{rema}
Up to unitary equivalence, the representation 
$\pi_n$ does not depend on the choice of a representative 
$n  \in F^{\times}/\varphi(F^{\times})$.
\end{rema}

\begin{proof}
Define
$$
\mathcal{H}_0 := \left\{ \phi \in \mathcal H \left| \phi ((x, 1)g) = \phi(g) \right. \right\},
$$
and let $\mathcal H_1$ be the orthogonal complement of $\mathcal H_0$. 
It is straightforward to prove that 
$$
\mathcal H_0\cong \mathsf L^2(T(F)\backslash T(\bA_F)).
$$ 
The following two lemmas prove that 
$$
\mathcal H_1 \cong \widehat{\bigoplus}_{n \in F^{\times}/\varphi(F^{\times})} \pi_n.
$$
\end{proof}

\begin{lemm} 
\label{proj}
For any $\phi \in \mathsf L^1(G(F)\backslash G(\bA_F))\cap \mathcal H$, the projection of 
$\phi$ onto $\mathcal H_0$ is given by
$$ 
\phi_0(g) := \int_{N(F) \backslash N(\bA_F)} \phi((x, 1)g)dx.
$$
\end{lemm}

\begin{proof}
It is easy to check that $\phi_0 \in \mathcal H_0$. Also, for any $\phi' \in \mathcal H_0$, we have
$$ 
\int_{G(F) \backslash G(\bA_F)} (\phi-\phi_0)\phi' \mathrm d g =0.
$$
\end{proof}

\begin{lemm}
\label{lemm:embedding}
We have
$$
\mathcal H_1 \cong \widehat{\bigoplus}_{n \in F^{\times}/\varphi(F^{\times})} \pi_{n}.
$$
\end{lemm}

\begin{proof}
First we note that the underlying Hilbert space of $\pi_{n}$ is $\mathsf L^2(W \backslash T(\bA_F))$, and that the group action is given by
$$ 
(x, a)\cdot f(b) = \psi_n(\varphi(b)x)f(ab),
$$
where $f$ is a square-integrable function on $T(\bA_F)$.
For $\phi \in C_c^\infty(G(F)\backslash G(\bA_F))\cap \mathcal H_1$, define
$$
f_{n, \phi}(a) := \int_{N(F)\backslash N(\bA_F)} 
\phi(x,a)\overline{\psi}_n(x) \, \mathrm  dx.
$$
Then, 
\begin{align*}
\parallel \phi \parallel_{\mathsf L^2}^2 &= \int_{T(F)\backslash T(\bA_F)}\int_{N(F)\backslash N(\bA_F)}\left|\phi(x, a)) \right|^2\mathrm dx\mathrm da^\times \\
&= \int_{T(F)\backslash T(\bA_F)} \sum_{\alpha \in F} \left|\int_{N(F)\backslash N(\bA_F)} \phi(x,a)\overline{\psi}(\alpha x)\mathrm  dx\right|^2 \mathrm da^\times \\
&= \int_{T(F)\backslash T(\bA_F)} \sum_{\alpha \in F^\times} \left|\int_{N(F)\backslash N(\bA_F)} \phi(x,a)\overline{\psi}(\alpha x)\mathrm dx\right|^2 \mathrm da^\times \\
&= \int_{T(F)\backslash T(\bA_F)} \sum_{\alpha \in F^\times} \sum_{n \in F^{\times}/\varphi(F^{\times})} \frac{1}{\# W}\left|\int_{N(F)\backslash N(\bA_F)} \phi(x,a)\overline{\psi}(n\varphi(\alpha) x)\mathrm dx\right|^2 \mathrm da^\times \\
&= \int_{T(F)\backslash T(\bA_F)} \sum_{\alpha \in F^\times} \sum_{n \in F^{\times}/\varphi(F^{\times})} \frac{1}{\# W}\left|\int_{N(F)\backslash N(\bA_F)} \phi(x, \alpha a)\overline{\psi}(nx)\mathrm dx\right|^2 \mathrm da^\times \\
&= \sum_{n \in F^{\times}/\varphi(F^{\times})} \frac{1}{\# W} \int_{T(\bA_F)} \left|\int_{N(F)\backslash N(\bA_F)} \phi(x, a)\overline{\psi}_n(x)dx\right|^2 \mathrm da^\times\\
& = \sum_{n \in F^{\times}/\varphi(F^{\times})} \parallel f_{n, \phi} \parallel_{\mathsf L^2}^2.
\end{align*}
The second equality is the Plancherel theorem for 
$N(F) \backslash N(\bA_F)$. Third equality follows from the previous lemma. 
The fourth equality follows from the left $G(F)$-invariance of $\phi$.
Thus, we obtain an unitary operator:
$$
I : \mathcal H_1 \rightarrow \widehat{\bigoplus}_{n \in F^{\times}/\varphi(F^{\times})} \pi_{n}.
$$
Compatibility with the group action is straightforward, 
so $I$ is actually a morphism of unitary representations.
We construct the inverse map of $I$ explicitly. 
For $f \in \mathsf C_c^{\infty}(W \backslash T(\bA_F))$, define
$$
\phi_{n, f}(x, a) := \frac{1}{\# W} \sum_{\alpha \in F^\times} \psi_n(\varphi(\alpha) x )f(\alpha a).
$$ 
The orthogonality of characters implies that 
$$
\begin{array}{cl}
{} & \int_{N(F)\backslash N(\bA_F)}\phi_{n, f}(x, a)\cdot \overline{\phi_{n, f}(x, a)}\,\, \mathrm dx \\
& \\
& = \int_{N(F)\backslash N(\bA_F)}(\sum_{\alpha \in F^\times}  \psi_n(\varphi(\alpha) x)f(\alpha a)) \cdot  
(\sum_{\alpha \in F^\times} \overline{\psi}_n(\varphi(\alpha) x)\overline{f}(\alpha a))\,\,\mathrm dx \\
& \\
& = \sum_{\alpha \in F^\times}|f(\alpha a)|^2.
\end{array}
$$
Substituting, we obtain
\begin{align*}
\parallel \phi_{n, f} \parallel^2 &= \int_{T(F)\backslash T(\bA_F)}\int_{N(F)\backslash N(\bA_F)}\left|\phi_{n, f}(x, a)\right|^2\mathrm dx\mathrm da^\times \\
&= \frac{1}{\# W}\int_{T(F)\backslash T(\bA_F)}\sum_{\alpha \in F^\times}|f(\alpha a)|^2\mathrm da^\times = \parallel f \parallel_n^2.
\end{align*}
Lemma \ref{proj} implies that that $\phi_f \in \mathcal H_1$ and we obtain a morphism 
$$
\Theta :  \widehat{\bigoplus}_{n \in F^{\times}/\varphi(F^{\times})} \pi_{n} \rightarrow \mathcal H_1.
$$ 
Now we only need to check that 
$\Theta I = \textit{id}$ and $I \Theta = \textit{id}$. The first follows from the Poisson formula: 
For any $\phi \in \mathsf C_c^\infty(G(F)\backslash G(\bA_F))\cap \mathcal H_1$,
\begin{align*}
\Theta I\phi &= \sum_{n \in F^{\times}/\varphi(F^{\times})} \frac{1}{\# W}\sum_{\alpha \in F^\times} \psi_n(\varphi(\alpha) x )\int_{N(F)\backslash N(\bA_F)} \phi(y, \alpha a)\overline{\psi}_n(y)\mathrm dy \\
&= \sum_{n \in F^{\times}/\varphi(F^{\times})} \frac{1}{\# W}\sum_{\alpha \in F^\times} \int_{N(F)\backslash N(\bA_F)} \phi(\varphi(\alpha)y, \alpha a) \overline{\psi}_n(\varphi(\alpha) (y-x))\mathrm dy \\
&= \sum_{n \in F^{\times}/\varphi(F^{\times})} \frac{1}{\# W}\sum_{\alpha \in F^\times} \int_{N(F)\backslash N(\bA_F)} \phi((y +x, a))\overline{\psi}_n(\varphi(\alpha) y)\mathrm dy \\
&= \sum_{\alpha \in F} \int_{N(F)\backslash N(\bA_F)} \phi((y, 1)(x, a))\overline{\psi}(\alpha y)\mathrm dy = \phi(x, a).
\end{align*}
The other identity, $I\Theta = \textit{id}$ is checked by a similar computation.
\end{proof}

To simplify notation, we now restrict to $F=\bQ$. For our applications in Sections~\ref{sect:plane} and \ref{sect:general}, 
we need to know an explicit orthonormal basis for the unique infinite-dimensional representation 
$\pi = \mathsf L^2(\bA_{\bQ}^{\times})$ of $G = G_1$. 
For any $n \geq 1$, define compact subgroups of $G(\bZ_p)$
$$
G(p^n \bZ_p) := \{ (x, a) \, |\, x \in p^n \bZ_p, \, a \in 1 + p^n \bZ_p \}.
$$
Let  $v_p : \bQ_p \ra \bZ$ be the discrete valuation on $\bQ_p$.

\begin{lemm}
\label{lemm:basis}
Let $\mathbf K_p = G(p^n \bZ_p)$. 
\begin{itemize}
\item When $n=0$, an orthonormal basis for 
$\mathsf L^2(\bQ_p^{\times})^{\mathbf K_p}$ is given by
$$
\{ \boldsymbol{1}_{p^j\bZ_p^{\times}} \, | \, j \geq 0 \}.
$$

\item When $n \geq 1$, an orthonormal basis for 
$\mathsf L^2(\bQ_p^{\times})^{\mathbf K_p}$ is given by
$$
\{ \lambda_p(\cdot / p^j)\boldsymbol{1}_{p^j\bZ_p^{\times}} \, | \, j \geq -n, \, \lambda_p \in \mathsf M_p \},
$$
where $\mathsf M_p$ is the set of multiplicative characters on $\bZ_p^{\times}/(1 + p^n\bZ_p)$.
\end{itemize}
Moreover, let $\mathbf K_{\textnormal{fin}} = \prod_p \mathbf K_p$ where $\mathbf K_p = G(p^{n_p} \bZ_p)$ and $n_p = 0$ for almost all $p$. Let $S$ be the set of primes with $n_p \neq 0$ and $N = \prod_p p^{n_p}$. Then an orthonormal basis for 
$\mathsf L^2(\bA_{\bQ, \textnormal{fin}}^{\times})^{\mathbf K_{\textnormal{fin}}}$ is given by
$$
\{ \otimes_{p\in S}\lambda_p(a_p \cdot p^{-v_p(a_p)}) \boldsymbol{1}_{\frac{m}{N}\bZ_p^{\times}}(a_p) \otimes'_{p \notin S} \boldsymbol{1}_{m\bZ_p^{\times}}(a_p)\, |\, m \in \bN ,\, \lambda_p \in \mathsf M_p\}.
$$
\end{lemm}

\begin{proof}
For the first assertion, let $f \in \mathsf L^2(\bQ_p^{\times})^{\mathbf K_p}$ where $\mathbf K_p = G(\bZ_p)$. Since it is $\mathbf K_p$-invariant, we have
$$
f(b_p \cdot a_p) = f(a_p),
$$
for any $b \in \bZ_p^{\times}$. Hence $f$ takes the form of
$$
f = \sum_{j = -\infty}^{\infty} c_j \boldsymbol{1}_{p^j\bZ_p^{\times}}
$$
where $c_j = f(p^j)$ and $\sum_{j = -\infty}^{\infty} |c_j|^2 < +\infty$.
On the other hand we have
$$
\psi_p(a_p\cdot x_p) f(a_p) = f(a_p)
$$
for any $x_p \in \bZ_p$. This implies that $f(p^j) = 0$ for any $j <0$. Thus the first assertion follows. The second assertion is treated similarly. The last assertion follows from the first and the second assertions.
\end{proof}

We denote these vectors by $\mathbf v_{m, \lambda}$ where $m \in \bN$ and $\lambda \in \mathsf M := \prod_{p \in S} \mathsf M_p$. Note that $\mathsf M$ is a finite set. Also we define
\begin{align*}
\theta_{m, \lambda, t}(g) &:= \Theta(\mathbf v_{m, \lambda} \otimes |\cdot|_{\infty}^{it})(g)\\
&= \sum_{\alpha \in \bQ^{\times}} \psi(\alpha x) \mathbf v_{m, \lambda}(\alpha a_{\textnormal{fin}}) \,  |\alpha a_\infty|_{\infty}^{it}.
\end{align*}

The following proposition is a combination of Lemma \ref{lemm:basis} and the standard Fourier analysis on the real line:
\begin{prop}
\label{prop:spectral}
Let $f \in \mathcal H_1^{\mathbf K}$. Suppose that
\begin{enumerate}
\item $I(f)$ is integrable, i.e.,
$$
I(f) \in \mathsf L^2(\bA^{\times})^{\mathbf K} \cap \mathsf L^1(\bA^{\times}),
$$
\item the Fourier transform of $f$ is also integrable i.e.
$$
\int_{-\infty}^{+\infty} |(f, \theta_{m, \lambda, t})| \, \mathrm dt < + \infty,
$$
for any $m \in \bN$ and $\lambda \in \mathbf M$.
\end{enumerate}
Then we have
$$
f(g) = \sum_{\lambda \in \mathbf M} \sum_{m=1}^{\infty} \frac{1}{4\pi} \int_{-\infty}^{+\infty} (f, \theta_{m, \lambda, t}) \theta_{m, \lambda, t}(g)\, \mathrm dt \, \, \textnormal{ a.e.,}
$$
where
$$
(f, \theta_{m, \lambda, t}) = \int_{G(\bQ) \backslash G(\bA_{\bQ})} f(g) \overline{\theta}_{m, \lambda, t}(g) \, \mathrm dg.
$$
\end{prop}

\begin{proof}
For simplicity, we assume that $n_p =0$ for all primes $p$. Let $I(f) = h \in \mathsf L^2(\bA^{\times})^{\mathbf K} \cap \mathsf L^1(\bA^{\times})$. It follows from the proof of Lemma \ref{lemm:embedding} that
$$
f(g) = \Theta(h)(g) = \sum_{\alpha \in \bQ^{\times}} \psi(\alpha x) h(\alpha a).
$$
Note that this infinite sum exists in both $\mathsf L^1$ and $\mathsf L^2$ sense. It is easy to check that
$$ 
\int_{\bA^{\times}} h(a) \mathbf v_{m}(a_{\text{fin}}) |a_\infty|_{\infty}^{-it}\, \mathrm da^{\times}
= (f, \theta_{m, t}).
$$
Write
$$
h = \sum_m \mathbf v_m \otimes h_m,
$$
where $h_m \in \mathsf L^2(\bR_{>0}, \mathrm da_\infty^{\times})$. The first and the second assumptions imply that $h_m$ and the Fourier transform of $\widehat{h}_m$ both are integrable.
Hence the inverse formula of Fourier transformation on the real line implies that
$$
h(a) = \sum_m \frac{1}{4\pi} \int_{-\infty}^{+\infty} (f, \theta_{m, t}) \mathbf v_{m}(a_{\text{fin}}) |a_\infty|_{\infty}^{it} \, \mathrm dt \, \, \, \text{ a.e.}.
$$
Apply $\Theta$ to both sides, and our assertion follows.
\end{proof}

We recall some results regarding Igusa integrals with rapidly 
oscillating phase, studied in \cite{CT-additive}:

\begin{prop}
\label{prop:osillatory}
Let $p$ be a finite place of $\bQ$  and $d, e \in \bZ$.
Let
$$
\Phi : \bQ_p^2 \times \bC^2 \ra \bC,
$$ 
be a function such that for each $(x,y)\in \bQ_p^2$, $\Phi((x,y), \mathbf s)$
is holomorphic in $\mathbf s = (s_1, s_2) \in \bC^2$.
Assume that the function $(x, y) \mapsto \Phi(x, y, \mathbf s)$ belongs to 
a bounded subset of the space of smooth compactly supported functions when 
$\Re (\mathbf s)$ belongs to a fixed compact subset of $\bR^2$.
Let $\Lambda$ be the interior of a closed convex cone generated by
$$
(1, 0), (0, 1), (d, e).
$$
Then, for any $\alpha \in \bQ_p^{\times}$, 
$$
\eta_\alpha (\mathbf s) = \int_{\bQ_p^2} |x|_p^{s_1} |y|_p^{s_2} \psi_p(\alpha x^d y^e) \Phi(x, y, \mathbf s) \mathrm dx_p^{\times} \mathrm dy_p^{\times},
$$
is holomorphic on $\mathsf T_\Lambda$. The same argument holds for the infinite place when $\Phi$ is a smooth function with compact supports.
\end{prop}

\begin{proof}
For the infinite place, use integration by parts and apply the convexity principle. For finite places, assume that $d, e$ are both negative. Let $\delta(x, y) = 1$ if $|x|_p =|y|_p =1$ and $0$ else. Then we have
\begin{align*}
\eta_\alpha (\mathbf s) &= \sum_{n, m \in \bZ} \int_{\bQ_p^2} |x|_p^{s_1} |y|_p^{s_2} \psi_p(\alpha x^d y^e) \Phi(x, y, \mathbf s) \delta (p^{-n}x, p^{-m}y) \, \mathrm dx_p^{\times} \mathrm dy_p^{\times}\\
&= \sum_{n, m \in \bZ} p^{-(ns_1 + ms_2)} \cdot \eta_{\alpha, n, m} (\mathbf s), 
\end{align*}
where
$$
\eta_{\alpha, n, m} (\mathbf s) = \int_{|x|_p = |y|_p =1}\!\!\!\!\!\!\!\!\!\! \!\!\!\!\!\!\!\!\!\!\psi_p(\alpha p^{nd + me} x^d y^e) \Phi(p^nx, p^my, \mathbf s) \, \mathrm dx_p^{\times} \mathrm dy_p^{\times}.
$$
Fix a compact subset of $\bC^2$ and assume that $\Re (\mathbf s)$ is in that compact set. The assumptions in our proposition mean that the support of $\Phi (\cdot, \mathbf s)$ is contained in a fixed compact set in $\bQ_p^2$, so there exists an integer $N_0$ such that $\eta_{\alpha, n, m}(\mathbf s) = 0$ if $n< N_0$ or $m < N_0$. Moreover our assumptions imply that there exists a positive real number $\delta$ such that $\Phi(\cdot, \mathbf s)$ is constant on any ball of radius $\delta$ in $\bQ_p^2$. This implies that if $1/p^n < \delta$, then for any $u \in \bZ_p^{\times}$,
\begin{align*}
\eta_{\alpha, n, m} (\mathbf s) &= \int \psi_p(\alpha p^{nd + me} x^d y^e u^d) \Phi(p^nxu, p^my, \mathbf s) \, \mathrm dx_p^{\times} \mathrm dy_p^{\times}\\
&=\int \psi_p(\alpha p^{nd + me} x^d y^e u^d) \Phi(p^nx, p^my, \mathbf s) \, \mathrm dx_p^{\times} \mathrm dy_p^{\times}\\
&= \int \int_{\bZ_p^{\times}} \psi_p(\alpha p^{nd + me} x^d y^e u^d) \, \mathrm du^{\times} \Phi(p^nx, p^my, \mathbf s) \, \mathrm dx_p^{\times} \mathrm dy_p^{\times},
\end{align*}
and the last integral is zero if $n$ is sufficiently large because of \cite[Lemma 2.3.5]{CT-additive}. Thus we conclude that there exists an integer $N_1$ such that $\eta_{\alpha, n, m}(\mathbf s) =0$ if $n >N_1$ or $m>N_1$. Hence we obtained that
$$
\eta_\alpha (\mathbf s) = \sum_{N_0 \leq n, m \leq N_1} p^{-(ns_1 + ms_2)} \cdot \eta_{\alpha, n, m} (\mathbf s),
$$
and this is holomorphic everywhere.

The case of $d < 0$ and $e=0$ is treated similarly.

Next assume that $d < 0$ and $e > 0$. Then again we have a constant $c$ such that $\eta_{\alpha, n, m} (\mathbf s) = 0$ if $1/p^n < \delta$ and $n|d| - me > c$. We may assume that $c$ is sufficiently large so that the first condition is unnecessary. Then we have
\begin{align*}
|\eta_\alpha (\mathbf s)| & \leq \sum_{N_0 \leq n} \sum_m p^{-\frac{n}{e}(e\Re (s_1) + |d|\Re (s_2))} \cdot p^{\frac{(n|d| - me)}{e}\Re (s_2)} \cdot | \eta_{\alpha, n, m} (\mathbf s) |\\
&\leq \sum_{N_0 \leq n} p^{-\frac{n}{e}(e\Re (s_1) + |d|\Re (s_2))} \cdot \frac{p^{\frac{c}{e}\Re (s_2)}}{1 - p^{-\frac{\Re(s_2)}{e}}}
\end{align*}
Thus $\eta_\alpha (\mathbf s)$ is holomorphic on $\mathsf T_\Lambda$.

\end{proof}
From the proof of Proposition \ref{prop:osillatory}, we can claim more for finite places:

\begin{prop}
\label{prop:osillatorybounds}
Let $\epsilon> 0$ be any small positive real number. Fix a compact subset $K$ of $\Lambda$, and assume that $\Re(\mathbf s)$ is in $K$. 
Define:
$$
\kappa (K) := 
\begin{cases}
\max \left\{ 0, -\frac{\Re(s_1)}{|d|}, -\frac{\Re(s_2)}{|e|} \right\} & \textnormal{ if $d < 0$ and $e < 0$,}\\
\max \left\{ 0, -\frac{\Re(s_1)}{|d|} \right\} & \textnormal{ if $d < 0$ and $e \geq 0$,}.
\end{cases}
$$
Then we have
$$
|\eta_\alpha (\mathbf s)| \ll
1/|\alpha|_p^{\kappa(K) + \epsilon}
$$
as $|\alpha|_p \ra 0$.
\end{prop}
\begin{proof}
Let $|\alpha|_p = p^{-k}$, and assume that both $d, e$ are negative. By changing variables, 
if necessary, we may assume that $N_0$ in the proof of Proposition \ref{prop:osillatory} is zero. If $k$ is sufficiently large, then one can prove that there exists a constant $c$ such that $\eta_{\alpha, n, m} (\mathbf s) = 0$ if $n|d| + m|e| > k+c$. Also it is easy to see that
$$
|p^{-(ns_1 + ms_2)}| \leq p^{(n|d| + m|e|)\kappa(K)}.
$$
Hence we can conclude that
$$
|\eta_\alpha (\mathbf s)| \ll k^2 1/|\alpha|_p^{\kappa(K)} \ll 1/|\alpha|_p^{\kappa(K) + \epsilon}.
$$

The case of $d < 0$ and $e = 0$ is treated similarly.

Assume that $d < 0$ and $e > 0$. Then we have a constant $c$ such that $\eta_{\alpha, n, m} (\mathbf s) = 0$ if $n|d| - me > k+c$. Thus we can conclude that
\begin{align*}
|\eta_\alpha (\mathbf s)| &\leq \sum_{m \geq 0} \sum_{n \geq 0} p^{-(n\Re (s_1) + m\Re (s_2))} |\eta_{\alpha, n, m} (\mathbf s)|\\
&\ll \sum_{m \geq 0} p^{-m\Re(s_2)} (me + k) p^{(me + k)\kappa (K, s_2)}\\
&\ll k 1/|\alpha|_p^{\kappa(K , s_2)} \sum_{m \geq 0} (m+1) p^{-m(\Re(s_2) - e\kappa (K, s_2)) }\\
&\ll 1/|\alpha|_p^{\kappa(K ,s_2) + \epsilon}.
\end{align*}
where 
$$
\kappa(K, s_2) = \max \left\{ 0, -\frac{\Re(s_1)}{|d|} : (\Re(s_1), \Re (s_2)) \in K \right\}.
$$
Thus we can conclude that
$$
|\eta_\alpha (\mathbf s)| \ll 1/|\alpha|_p^{\kappa(K ,s_2) + \epsilon} \ll 1/|\alpha|_p^{\kappa(K) + \epsilon}.
$$
\end{proof}

\section{The projective plane}
\label{sect:plane}

In this section, we implement the program described in Section~\ref{sect:height-zeta} for 
the simplest equivariant compactifications of $G=G_1 = \mathbb G_a\rtimes \mathbb G_m$, namely, 
the projective plane $\bP^2$, for a {\em one-sided}, right, action of $G$ given by  
$$
G \ni (x, a) \mapsto [x_0:x_1:x_2]=(a:a^{-1}x:1) \in \bP^2.
$$ 
The boundary consists of two lines, $D_0$ and $D_2$ given by the 
vanishing of  $x_0$ and $x_2$. We will use the following identities:
\begin{align*}
&\text{div}(a) = D_0 - D_2,\\
&\text{div}(x) = D_0+D_1 -2D_2,\\
&\text{div}(\omega) = -3D_2,
\end{align*}
where $D_1$ is given by the vanishing of $x_1$ 
and $\omega$ is the right invariant top degree form.
The height functions are given by
\begin{align*}
&\sH_{D_0, p}(a,x) = \frac{\max \{ |a|_p, |a^{-1}x|_p, 1\}}{|a|_p}, &\sH_{D_2, p}(a,x)=\max \{|a|_p, |a^{-1}x|_p, 1\},\\
&\sH_{D_0, \infty}(a,x) = \frac{\sqrt{|a|^2 + |a^{-1}x|^2 + 1}}{|a|},&\sH_{D_2, p}(a,x)=\sqrt{|a|^2+|a^{-1}x|^2+1},\\
&\sH_{D_0} = \prod_p \sH_{D_0, p} \times \sH_{D_0, \infty}, &\sH_{D_2} = \prod_p \sH_{D_2, p} \times \sH_{D_2, \infty},
\end{align*}
and the height pairing by
$$
\sH(\mathbf{s}, g) = \sH_{D_0}^{s_{0}}(g) \sH_{D_2}^{s_{2}}(g),
$$
for $\mathbf{s} = s_{0} D_0 + s_{2}D_2$ and $g \in G(\bA)$.
The height zeta function takes the form
$$
\mathsf Z(\mathbf{s}, g) = \sum_{\gamma \in G(\bQ)} \sH(\mathbf{s}, \gamma g)^{-1}.
$$
The proof of Northcott's theorem shows that the 
Dirichlet series $\mathsf Z(\mathbf{s}, g)$ 
converges absolutely and normally to a holomorphic function, for 
$\Re(\mathbf{s})$ is sufficiently large,  which is 
continuous in $g \in G(\bA)$. Moreover, if 
$\Re(\mathbf{s})$ is sufficiently large, then 
$$
\sZ(\mathbf{s}, g) \in \mathsf L^2(G(\bQ)\backslash G(\bA)) \cap 
\mathsf L^1(G(\bQ)\backslash G(\bA)).
$$
According to Proposition \ref{prop:parameter}, 
we have the following decomposition:
$$
\mathsf L^2(G(\bQ)\backslash G(\bA)) = \mathsf L^2(\mathbb G_m(\bQ)\backslash \mathbb G_m(\bA)) \oplus \pi,
$$
and we can write
$$
\sZ(\mathbf{s}, g) = \sZ_0(\mathbf{s}, g) + \sZ_1(\mathbf{s}, g).
$$
The analysis of $\sZ_0(\mathbf{s}, \text{id})$ is a special case of 
our considerations in Section~\ref{sect:height-zeta}, in particular Theorem~\ref{theo:mainterm}
(for further details, see \cite{BT98} and \cite{volume}).
The conclusion here is that there exist  a $\delta > 0$ and 
a function $\mathsf h$ which is holomorphic on the tube domain $\mathsf{T}_{>3-\delta}$ such that
$$
\sZ_0(\mathbf{s}, \text{id}) = \frac{\mathsf h(s_{0} +s_{2})}{(s_{0} +s_{2} -3)}.
$$
The analysis of $\sZ_1(\mathbf{s}, \text{id})$, i.e., of the contribution 
from the unique infinite-dimensional representation occurring in 
$\mathsf L^2(G(\bQ)\backslash G(\bA))$,
is the main part of this section. 
Define
$$
\mathbf K = \prod_p \mathbf K_p \cdot \mathbf K_{\infty} = \prod_p G(\bZ_p) \cdot \{ (0, \pm 1)\}.
$$
Since the height functions are $\bK$-invariant, 
$$
\sZ_1(\mathbf{s}, g) \in \pi^{\bK} \simeq \mathsf L^2 (\bA^{\times})^{\bK}.
$$
Lemma~\ref{lemm:basis} provides a choice of an orthonormal basis for 
$\mathsf L^2 (\bA^{\times}_{{\rm fin}})$. Combining with the Fourier expansion at the archimedean place, 
we obtain the following spectral expansion of $\sZ_1$:

\begin{lemm}
\label{lemm:transformation}
Assume that $\Re(\mathbf{s})$ is sufficiently large. Then 
$$
\sZ_1(\mathbf{s}, g) = 
\sum_{m \geq 1} \frac{1}{4\pi}\int_{-\infty}^{\infty} (\sZ(\mathbf{s}, g), 
\theta_{m,t}(g))\theta_{m,t}(g)\mathrm dt,
$$
where $\theta_{m, t}(g) = \Theta(\mathbf v_m \otimes |\cdot|^{it})(g)$.
\end{lemm}
\begin{proof}
See Lemma \ref{lemm:generalspectral}.
\end{proof}

It is easy to see that
\begin{align*}
(\sZ(\mathbf{s}, g), \theta_{m,t}(g)) &= \int_{G(\bQ)\backslash G(\bA)} \sZ(\mathbf{s}, g)\theta_{m,t}(g)\mathrm d g\\
&= \int_{G(\bA)} \sH(\boldsymbol{s}, g)^{-1}\bar{\theta}_{m,t}(g)\mathrm d g\\
&= \sum_{\alpha \in \bQ^{\times}} \int_{G(\bA)} \sH(\mathbf{s}, g)^{-1} \bar{\psi}(\alpha x) \mathbf v_m (\alpha a_{\text{fin}}) |\alpha a_{\infty}|^{-it} \mathrm d g\\
&= \sum_{\alpha \in \bQ^{\times}} \prod_p \sH_p'(\mathbf{s}, m, \alpha) \cdot 
\sH_{\infty}'(\mathbf{s}, t, \alpha),
\end{align*}
where
\begin{align*}
&\sH_p'(\mathbf{s}, m, \alpha) = \int_{G(\bQ_p)}\sH_p(\mathbf{s}, g_p)^{-1}\bar{\psi}_p(\alpha x_p) \boldsymbol{1}_{m\bZ_p^{\times}} (\alpha a_p) \mathrm d g_p,\\
&\sH_{\infty}'(\mathbf{s}, t, \alpha) = \int_{G(\bR)}\sH_\infty(\mathbf{s}, g_\infty)^{-1}\bar{\psi}_\infty(\alpha x_\infty)  |\alpha a_{\infty}|^{-it} \mathrm d g_\infty.
\end{align*}
Note that $\theta_{m, t}(\text{id}) = 2|m|^{it}$. Hence we can conclude that
\begin{align*}
&\sZ_1(\mathbf{s}, \text{id}) = 
\sum_{\alpha \in \bQ^{\times}} \sum_{m=1}^\infty \frac{1}{2\pi}\int_{-\infty}^{\infty} \prod_p \sH_p'(\mathbf{s}, m, \alpha) \cdot 
\sH_{\infty}'(\mathbf{s}, t, \alpha) |m|^{it} \mathrm dt\\
& \\
&= \sum_{\alpha \in \bQ^{\times}} \!\sum_{m=1}^\infty \! \frac{1}{2\pi}\int_{-\infty}^{\infty}\!\! \prod_p \!\int_{G(\bQ_p)}\!\!\!\!\!\sH_p(\mathbf{s}, g_p)^{-1}\bar{\psi}_p(\alpha x_p) \boldsymbol{1}_{m\bZ_p^{\times}} (\alpha a_p) |\alpha a|_p^{-it}\mathrm d g_p \cdot \sH_{\infty}'(\mathbf{s}, t, \alpha) \mathrm dt\\
& \\
&= \sum_{\alpha \in \bQ^{\times}} \frac{1}{2\pi}\int_{-\infty}^{\infty} \prod_p \widehat{\sH}_p(\mathbf{s}, \alpha, t) \cdot 
\widehat{\sH}_\infty (\mathbf{s}, \alpha, t) \mathrm dt,
\end{align*}
where
\begin{align*}
&\widehat{\sH}_p(\mathbf{s}, \alpha, t) = \int_{G(\bQ_p)}\sH_p(\mathbf{s}, g_p)^{-1}\bar{\psi}_p(\alpha x_p) \boldsymbol{1}_{\bZ_p} (\alpha a_p) |a_p|_p^{-it}\mathrm d g_p,\\
&\widehat{\sH}_\infty (\mathbf{s}, \alpha, t)=\int_{G(\bR)}\sH_\infty(\mathbf{s}, g_\infty)^{-1}\bar{\psi}_\infty(\alpha x_\infty)  |a_{\infty}|^{-it} \mathrm d g_\infty.
\end{align*}
It is clear that 
$$
\widehat{\sH}_p(\mathbf{s}, \alpha, t) = \widehat{\sH}_p((s_{0} - it, s_{2} +it),  \alpha, 0),
$$ 
so we only need to study 
$\widehat{\sH}_p(\mathbf{s}, \alpha) = \widehat{\sH}_p(\mathbf{s}, \alpha, 0)$.  
To do this, we introduce some notation. 
We have the canonical integral model of $\bP^2$ over $\Spec(\bZ)$, 
and for any prime $p$, we have the reduction map modulo $p$:
$$
\rho :G(\bQ_p) \subset \bP^2(\bQ_p) = \bP^2(\bZ_p) \rightarrow \bP^2(\bF_p)
$$
This is a continuous map from $G(\bQ_p)$ to $\bP^2(\bF_p)$. Consider the following open sets:
\begin{align*}
U_{\phi}  & = \rho^{-1}(\bP^2 \setminus (D_0\cup D_2))= \{ |a|_p =1, |x|_p \leq 1 \}\\
U_{D_0}   & = \rho^{-1}(D_0 \setminus (D_0 \cap D_2))   = \{|a|_p < 1, |a^{-1}x|_p \leq 1 \}\\
U_{D_2}  & = \rho^{-1}(D_2 \setminus (D_0 \cap D_2)) = \{|a|_p >1, |a^{-2}x|_p  \leq 1\}\\
U_{D_0, D_2} &  = \rho^{-1}(D_0 \cap D_2)              = \{|a^{-1}x|_p >1,  |a^{-2}x|_p > 1\}.
\end{align*}
The height functions have a partial left invariance, i.e., 
they are invariant under the left action of the 
compact subgroup $\{ (0, b) \mid b \in \bZ_p^{\times} \}$. This implies that
$$
\widehat{\sH}_p(\mathbf{s}, \alpha) =  \int_{G(\bQ_p)}\sH_p(\mathbf{s}, g)^{-1}\int_{\bZ_p^{\times}} \bar{\psi}_p(\alpha bx)\mathrm db^{\times}  \boldsymbol{1}_{\bZ_p} (\alpha a) \mathrm d g.
$$
We record the following useful lemma (see, e.g., \cite[Lemma 10.3]{VecIII}):
 
\begin{lemm}
\label{lemm:useful}
$$
\int_{\bZ_p^{\times}} \bar{\psi}_p(bx)\mathrm db^{\times} =
\begin{cases}
1 & \textnormal{ if } |x|_p \leq 1,\\
-\frac{1}{p-1} & \textnormal{ if } |x|_p =p,\\
0 & \textnormal{ otherwise}.
\end{cases}
$$
\end{lemm}

\begin{lemm}
Assume that $|\alpha|_p = 1$. Then 
$$
\widehat{\sH}_p(\mathbf{s}, \alpha) = 
\frac{\zeta_p(s_{0} +1)\zeta_p(2s_{0} + s_{2})}{\zeta_p(s_{0} + s_{2})}.
$$
\end{lemm}

\begin{proof}
We apply Lemma~\ref{lemm:useful} and obtain
\begin{align*}
\widehat{\sH}_p(\boldsymbol{s}, \alpha) &= \int_{U_{\emptyset}} + \int_{U_{D_0}} + \int_{U_{D_0, D_2}}\\
& \\
&= 1 + \frac{p^{-(s_{0}+1)}}{1-p^{-(s_{0} +1)}} + 
\frac{p^{-(2s_{0} + s_{2})} - p^{-(s_{0} + s_{2})}}{(1-p^{-(s_{0} +1)})(1-p^{-(2s_{0} +s_{2})})}\\
&  \\
&= \frac{\zeta_p(s_{0} +1)\zeta_p(2s_{0} + s_{2})}{\zeta_p(s_{0} + s_{2})}.
\end{align*}
\end{proof}

\begin{lemm}
Assume that $|\alpha|_p >1$. Let $|\alpha|_p = p^k$. Then 
$$
\widehat{\sH}_p(\mathbf{s}, \alpha) = p^{-k(s_{0}+1)} \widehat{\sH}_p(\mathbf{s}, 1).
$$
\end{lemm}
\begin{proof}
Use the lemma, and we obtain that
\begin{align*}
\widehat{\sH}_p(\mathbf{s}, \alpha) &= \int_{U_{D_0}} + \int_{U_{D_0, D_2}}\\
&= \frac{p^{-k(s_{0}+1)}}{1-p^{-(s_{0} +1)}} + p^{-k(s_{0}+1)} 
\frac{p^{-(2s_{0} + s_{2})} - p^{-(s_{0} + s_{2})}}{(1-p^{-(s_{0} +1)})(1-p^{-(2s_{0} +s_{2})})}\\
&= p^{-k(s_{0}+1)} \widehat{\sH}_p(\mathbf{s}, 1).
\end{align*}
\end{proof}

\begin{lemm}
Assume that $|\alpha|_p <1$. Let $|\alpha|_p = p^{-k}$. Then $\widehat{\sH}_p(\mathbf{s}, \alpha)$ 
is holomorphic on the tube domain $\mathsf T_\Lambda$ over the cone 
$$
\Lambda = \{s_{0} > -1, s_{0}+s_{2} > 0, 2s_{0} + s_{2} > 0 \}.
$$ 
Moreover, for any compact subset of $\mathsf \Lambda$, 
there exists a constant $C > 0$  such that
$$
|\widehat{\sH}_p(\mathbf{s}, \alpha)| \leq C k \max \{ 1, p^{-\frac{k}{2}\Re(s_{2} -2)} \}
$$
for any $\mathbf{s}$ with real part in this compact set. 
\end{lemm}

\begin{proof}
It is easy to see that
\begin{align*}
\widehat{\sH}_p(\mathbf{s}, \alpha) &= \int_{U_{\emptyset}} + \int_{U_{D_0}} + \int_{U_{D_2}} + \int_{U_{D_0, D_2}}\\
& \\
&= 1 + \frac{p^{-(s_{0}+1)}}{1-p^{-(s_{0} +1)}} + \int_{U_{D_2}} + \int_{U_{D_0, D_2}}.
\end{align*}
On $U_{D_2}$, we choose $(1:x_1:x_2)$ as coordinates, then we have
\begin{align*}
\int_{U_{D_2}} &= \int_{|x_1|_p \leq 1, |x_2|_p < 1} |x_2|_p^{s_{2} -2} \int_{\bZ_p^{\times}} 
\bar{\psi}_p(\alpha b \frac{x_1}{x_2^2})\mathrm db^{\times} \boldsymbol{1}_{\bZ_p} (\alpha x_2^{-1}) \mathrm dx_1 \,\mathrm dx_2^{\times}\\
&= \int_{p^{-\frac{k}{2}} \leq |x_2|_p <1} |x_2|_p^{s_{2} -2} \mathrm dx_2^{\times}
\end{align*}
Hence this integral is holomorphic everywhere, and we have
$$
\left|\int_{U_{D_2}} \right| < k \max \{1, p^{-\frac{k}{2}\text{Re}(s_{2}-2)} \}.
$$
On $U_{D_0, D_2}$, we choose $(x_0:1:x_2)$ as coordinates and obtain
\begin{align*}
& \frac{1}{1-p^{-1}}  \int_{U_{D_0, D_2}} = \int_{|x_0|_p,|x_2|_p <1}  |x_0|_p^{s_{0} +1} |x_2|_p^{s_{2} -2} \int_{\bZ_p^{\times}} \bar{\psi}_p(\alpha b \frac{x_0}{x_2^2})\mathrm db^{\times} \mathrm dx_0^{\times}\mathrm dx_2^{\times}\\
&= -\frac{p^{-1}}{1-p^{-1}}p^{(k+1)(s_{0} +1)}
\frac{p^{-(1-\left[-\frac{k}{2} \right] )(2s_{0} + s_{2})}}{1-p^{-(2s_{0} + s_{2})}} + 
\int |x_0|_p^{s_{0} +1} |x_2|_p^{s_{2}-2} \mathrm dx_0^{\times}\mathrm dx_2^{\times},
\end{align*}
where the last integral is
over 
$$
\{|x_0|_p <1, |x_2|_p <1, |\alpha x_0|_p \leq |x_2|_p^2\},
$$ 
and is equal to 
\begin{align*}
&\!\!\!\! \sum_{j, l \geq 1, 2j -l \leq k}p^{-l(s_{0} +1) -j(s_{2}-2)}\\
& \\
&= \sum_{j \leq \left[ \frac{k}{2} \right]}p^{-j(s_{2}-2)} \frac{p^{-(s_{0} +1)}}{1-p^{-(s_{0} +1)}} + \sum_{j > \left[ \frac{k}{2} \right]}p^{-j(s_{2}-2)} \frac{p^{-(2j -k)(s_{0} +1)}}{1-p^{-(s_{D_0} +1)}}\\
&  \\
&=\sum_{j \leq \left[ \frac{k}{2} \right]}p^{-j(s_{2}-2)} \frac{p^{-(s_{0} +1)}}{1-p^{-(s_{0} +1)}} +\frac{p^{-(\left[ \frac{k}{2} \right] +1)(2s_{0} + s_{2})}}{1-p^{-(2s_{0} + s_{2})}}\frac{p^{k(s_{0} +1)}}{1-p^{-(s_{0} +1)}}.
\end{align*}
From this we can see that 
$\int_{U_{D_0, D_2}}$ is holomorphic on $\mathsf T_{\Lambda}$
and that for any compact subset of $\Lambda$, 
we can find a constant $C>0$  such that
$$
\left|\int_{U_{D_0, D_2}}\right| < Ck \max \{ 1, p^{-\frac{k}{2}\Re(s_{2}-2)} \},
$$
for $\Re(\mathbf s)$ in this compact.
\end{proof}

Next, we study the local integral at the real place. Again, 
$$
\widehat{\sH}_\infty(\mathbf{s}, \alpha, t) = \widehat{\sH}_\infty((s_{0} - it, s_{2} +it), \alpha, 0),
$$ 
and we start with 
$$
\widehat{\sH}_\infty(\mathbf{s}, \alpha) = \widehat{\sH}_\infty(\mathbf{s}, \alpha, 0).
$$
\begin{lemm}
The function 
$$
\mathbf{s}\mapsto \widehat{\sH}_\infty(\mathbf{s}, \alpha),
$$
is holomorphic on the tube domain 
$\mathsf T_{\Lambda'}$ over 
$$
\Lambda' = \{ s_{0} > -1, s_{2} > 0, 2s_{0} +s_{2}>0\}.
$$ 
Moreover, for any $r\in \bN$ and any compact subset of 
$$
\Lambda_r' = \{ s_{0} > -1 + r, s_{2} > 0\},
$$ 
there exists a constant $C>0$ such that
$$
|\widehat{\sH}_\infty(\mathbf{s}, \alpha)| < \frac{C}{|\alpha|_\infty^r},
$$
for any $\mathbf{s}$ in the tube domain over this compact.
\end{lemm}

\begin{proof}
Let $U_{\emptyset} = X(\bR) \setminus (D_0 \cup D_2)$, $U_{D_i}$ be a small tubular neighborhood of 
$D_i$ minus $D_0\cap D_2$, and $U_{D_0, D_2}$ be a small neighborhood of $D_0 \cap D_2$. 
Then $\{U_{\emptyset}, U_{D_0}, U_{D_2}, U_{D_0, D_2}\}$ is an open covering of 
$X(\bR)$, and consider the partition of unity for this covering; 
$\theta_{\emptyset}, \theta_{D_0}, \theta_{D_2}, \theta_{D_0, D_2}$. 
Then we have
\begin{align*}
\widehat{\sH}_\infty(\mathbf{s}, \alpha) &= \int_{U_{\emptyset}}\sH_\infty^{-1}\bar{\psi}_\infty(\alpha x_\infty) \theta_{\emptyset} \mathrm d g_\infty + \int_{U_{D_0}} + \int_{U_{D_2}} + \int_{U_{D_0, D_2}}.
\end{align*}
On $U_{D_0, D_2}$, we choose $(x_0:1:x_2)$ as analytic coordinates and obtain
$$
\int_{U_{D_0, D_2}} = \int_{\bR^2} |x_0|^{s_{0} +1} |x_2|^{s_{2}-2} 
\bar{\psi}(\alpha \frac{x_0}{x_2^2}) 
\Phi (\boldsymbol{s}, x_0, x_2)\mathrm  dx_0^{\times}dx_2^{\times},
$$
where $\Phi$ is a smooth bounded function with compact support. 
Such oscillatory integrals have been studied in \cite{CT-additive}, in our case 
the integral is holomorphic if $\text{Re}(s_{0}) > -1$ and $\text{Re}(2s_{0} + s_{2}) > 0$. 
Assume that $\Re(\mathbf{s})$ is sufficiently large. Integration by parts implies that 
$$
\int_{U_{D_0, D_2}} = \frac{1}{\alpha^r} \int_{\bR^2} |x_0|^{s_{0} +1 - r} 
|x_2|^{s_{2}-2 + 2r} \bar{\psi}(\alpha \frac{x_0}{x_2^2}) 
\Phi' (\boldsymbol{s}, x_0, x_2) dx_0^{\times}\mathrm dx_2^{\times},
$$
and this integral is holomorphic if 
$\Re(s_{0}) > -1 + r$ and 
$\Re(s_{2}) > 2 - 2r$.
Thus, our second assertion follows. 
The other integrals are studied similarly.
\end{proof}

\begin{lemm}
For any compact set $K \subset \Lambda_2'$, there exists a constant $C>0$ such that
$$
|\widehat{\sH}_\infty(\mathbf{s}, \alpha, t)| < \frac{C}{|\alpha|^2(1+t^2)},
$$
for any $\mathbf{s} \in \mathsf T_K$.
\end{lemm}

\begin{proof}
Consider a left invariant differential operator $\partial_a = a\partial / \partial a$. 
Integrating by parts we obtain that
$$
\widehat{\sH}_\infty(\mathbf{s}, \alpha, t) = -\frac{1}{t^2} \int_{G(\bR)}\partial_a^2 
\sH_\infty(\mathbf{s}, g_\infty)^{-1}\bar{\psi}_\infty(\alpha x_\infty)  |a_{\infty}|^{-it} \mathrm d g_\infty,
$$
According to \cite{VecIII}, 
$$
\partial_a^2 \sH_\infty(\mathbf{s}, g_\infty)^{-1} = 
\sH_\infty(\mathbf{s}, g_\infty)^{-1} \times (\text{a bounded smooth function}),
$$ 
so we can apply the discussion of the previous proposition.
\end{proof}

\begin{lemm}
The Euler product
$$
\prod_p \widehat{\sH}_p(\mathbf{s}, \alpha, t) \cdot \widehat{\sH}_\infty (\mathbf{s}, \alpha, t),
$$
is holomorphic on the tube domain $\mathsf T_\Omega$ over 
$$
\Omega = \{s_{0} > 0, s_{2} >0, 2s_{0}+ s_{2} >1\}.
$$ 
Moreover, let $\alpha = \frac{\beta}{\gamma}$, where $\gcd(\beta, \gamma)=1$. 
Then for any $\epsilon > 0$ and any compact set 
$$
K \subset \Omega' = \{s_{0} > 1, s_{2} >0, 2s_{0}+ s_{2} >1\},
$$ 
there exists a constant $C>0$ such that
$$
|\prod_p \widehat{\sH}_p(\mathbf{s}, \alpha, t) \cdot \widehat{\sH}_\infty (\mathbf{s}, \alpha, t)| < C\cdot \frac{\max \{ 1, \sqrt{|\beta|}^{-\Re(s_{2} -2)} \}}{|\beta|^{2 - \epsilon}|\gamma|^{\Re(s_{0} -1)}},
$$
for all $\mathbf{s} \in \mathsf T_K$.
\end{lemm}

\begin{theo}
There exists  $\delta >0$ such that $\mathsf Z_1(s_{0} + s_{2}, \textnormal{id})$ is holomorphic on $\mathsf T_{> 3 -\delta}$.
\end{theo}
\begin{proof}
Let $\delta >0$ be a sufficiently small real number, and define 
$$
\Lambda = \{ s_{0} > 2 + \delta, s_{2} > 1 - 2\delta\}.
$$
It follows from the previous proposition that for any $\epsilon >0$ and any compact set $K \subset \Lambda$, there exists a constant $C>0$ such that
$$
|\prod_p \widehat{\sH}_p(\mathbf{s}, \alpha, t) \cdot \widehat{\sH}_\infty (\mathbf{s}, \alpha, t)| < \frac{C}{(1+t^2)|\beta|^{\frac{3}{2} -\epsilon -\delta}|\gamma|^{1+\delta}}.
$$
From this inequality, we can conclude that the integral
$$
\int_{-\infty}^{\infty} \prod_p \widehat{\sH}_p(\mathbf{s}, \alpha, t) \cdot \widehat{\sH}_\infty (\mathbf{s}, \alpha, t) \mathrm dt,
$$
converges uniformly and absolutely to a holomorphic function on $\mathsf T_K$. Furthermore, we have
$$
\left|\int_{-\infty}^{\infty} \prod_p \widehat{\sH}_p(\mathbf{s}, \alpha, t) \cdot \widehat{\sH}_\infty (\mathbf{s}, \alpha, t) \mathrm dt \right| < \frac{C'}{|b|^{\frac{3}{2} -\epsilon -\delta}|c|^{1+\delta}}.
$$
For sufficiently small $\epsilon>0$ and $\delta>0$, the sum
$$
\sum_{\alpha \in \bQ^{\times}} \frac{1}{2\pi}\int_{-\infty}^{\infty} \prod_p \widehat{\sH}_p(\mathbf{s}, \alpha, t) \cdot \widehat{\sH}_\infty (\mathbf{s}, \alpha, t) \mathrm dt,
$$
converges absolutely and uniformly to a function in $s_0+s_2$. 
This concludes the proof of our theorem.
\end{proof}

\section{Geometrization}
\label{sect:general}

In this section we geometrize the 
method described in Section~\ref{sect:plane}.
Our main theorem is:

\begin{theo}
\label{theo:general}
Let $X$ be a smooth projective equivariant compactification of 
$G=G_1$ over $\bQ$, under the right action. Assume that
the boundary divisor has strict normal crossings. 
Let $a, x \in \bQ(X)$ be rational functions, 
where $(x,a)$ are the standard coordinates on  
$G \subset X$. Let $E$ be the Zariski closure of $\{x=0\} \subset G$.
Assume that:
\begin{itemize}
\item the union of the boundary and $E$ is a divisor with strict normal crossings,
\item $\textnormal{div}(a)$ is a reduced divisor, and
\item for any pole $D_\iota$ of $a$, one has
$$
-\textnormal{ord}_{D_\iota} (x) > 1. 
$$
\end{itemize}
Then Manin's conjecture holds for $X$.
\end{theo}

The remainder of this section is devoted to a proof of this fact. 
Blowing up  the zero-dimensional subscheme 
$$
\text{Supp}(\text{div}_0(a)) \cap \text{Supp}(\text{div}_\infty(a)),
$$
if necessary, we may assume that 
$$
\text{Supp}(\text{div}_0(a)) \cap \text{Supp}(\text{div}_\infty(a)) =  \emptyset.
$$
The height functions 
are invariant under the right action of some compact subgroup 
$\mathbf K_p \subset G(\bZ_p)$. Moreover, we can assume that 
$\bK_p = G(p^{n_p}\bZ_p)$, for some $n_p\in \bZ_{\ge 0}$.
Let $S$ be the set of bad places for $X$; 
note that $n_p = 0$ for all $p \notin S$. For simplicity, 
we assume that the height function at the infinite place is 
invariant under the action of 
$\mathbf K_\infty = \{ (0, \pm1)\}$.

\begin{lemm}
\label{lemm:L2}
We have
$$
\sZ(\mathbf s, g) \in \mathsf L^2(G(\bQ)\backslash G(\bA))^{\mathbf K} \cap \mathsf L^1(G(\bQ)\backslash G(\bA)).
$$
\end{lemm}
\begin{proof}
First it is easy to see that
\begin{align*}
\int_{G(\bQ)\backslash G(\bA)} |\sZ (\mathbf s, g)| \, \mathrm dg
&\leq \int_{G(\bQ)\backslash G(\bA)} \sum_{\gamma \in G(\bQ)} |\sH (\mathbf s, \gamma g)|^{-1} \, \mathrm dg\\
&= \int_{G(\bA)} \sH ( \Re (\mathbf s), g)^{-1} \, \mathrm dg,
\end{align*}
and the last integral is bounded when $\Re(\mathbf s)$ is sufficiently large. (See \cite[Proposition 4.3.4]{volume}.) Hence it follows that $\sZ(\mathbf s, g)$ is integrable.
To conclude that $\sZ(\mathbf s, g)$ is square-integrable, we prove that $\sZ (\mathbf s, g) \in \mathsf L^{\infty}$ for $\Re (\mathbf s)$ sufficiently large. Let $u, v$ be sufficiently large positive real numbers. Assume that $\Re(\mathbf s)$ is in a fixed compact subset of $\Pic^G(X) \otimes \bR$ and sufficiently large. Then we have
$$
\sH ( \Re (\mathbf s), g)^{-1} \ll \sH_1(a)^{-u} \cdot \sH_2(x)^{-v}
$$
where
\begin{align*}
&\sH_1(a) = \prod_p \max \{ |a_p|_p, |a_p|_p^{-1} \} \cdot \sqrt{|a_\infty|_\infty^2 +|a_\infty|_\infty^{-2}}\\
&\sH_2(x) = \prod_p \max \{ 1, |x_p|_p \} \cdot \sqrt{1+|x_\infty|_\infty^2}.
\end{align*}
Since $\sZ(\mathbf s, g)$ is $G(\bQ)$-periodic, we may assume that $|a_p|_p=1$ where $g_p = (x_p, a_p)$.
Then we obtain that
\begin{align*}
\sum_{\gamma \in G(\bQ)} \sH ( \Re (\mathbf s), \gamma g)^{-1}
& \ll \sum_{\alpha \in \bQ^{\times}} \sum_{\beta \in \bQ} \sH_1(\alpha a)^{-u} \cdot \sH_2(\alpha x + \beta)^{-v}\\
&\leq \sum_{\alpha \in \bQ^{\times}} \sH_{1, \text{fin}}(\alpha)^{-u} \sZ_2(\alpha x),
\end{align*}
where
$$
\sZ_2(x) = \sum_{\beta \in \bQ} \sH_2(x + \beta)^{-v}.
$$
It is known that $\sZ_2$ is a bounded function for sufficiently large $v$, (see \cite{VecIII}) so we can conclude that $\sZ(\mathbf s, g)$ is also a bounded function because
$$
\sum_{\alpha \in \bQ^{\times}} \sH_{1, \text{fin}}(\alpha)^{-u} < +\infty,
$$
for sufficiently large $u$.
\end{proof}

By Proposition~\ref{prop:parameter}, 
the height zeta function decomposes as
$$
\sZ(\mathbf s, \text{id}) = 
\sZ_0(\mathbf s, \text{id}) + \sZ_1(\mathbf s, \text{id}).
$$
Analytic properties of 
$\sZ_0(\mathbf s, \text{id})$ were
established in Section~\ref{sect:height-zeta}. 
It remains to show 
that $\sZ_1(\mathbf s, \text{id})$ 
is holomorphic on a tube domain over an open neighborhood of the 
shifted effective cone $-K_X + \Lambda_{\text{eff}}(X)$. To conclude this, we use the spectral decomposition of $\sZ_1$:

\begin{lemm}
\label{lemm:generalspectral}
We have
$$
\sZ_1(\mathbf s, \textnormal{id}) = \sum_{\lambda \in \mathbf M} \sum_{m=1}^{\infty} \frac{1}{4\pi} \int_{-\infty}^{+\infty} 
(\sZ(\mathbf s, g), \theta_{m, \lambda, t}) \theta_{m, \lambda, t}(\textnormal{id})\, \mathrm dt.
$$
\end{lemm}
\begin{proof}
To apply Proposition \ref{prop:spectral}, we need to check that $\sZ_1$ satisfies the assumptions of Proposition \ref{prop:spectral}. The proof of Lemma \ref{lemm:embedding} implies that
$$
I(\sZ_1) = \int_{N(\bQ)\backslash N(\bA)} \sZ (\mathbf s, g) \psi(x) \, \mathrm dx.
$$
Thus we have 
\begin{align*}
\int_{T(\bA)}|I(\sZ_1)| \, \mathrm da^{\times} 
&= \int_{T(\bA)}\left| \int_{N(\bQ)\backslash N(\bA)} \sZ (\mathbf s, g) \psi(x) \, \mathrm dx \right| \, \mathrm da^{\times}\\
&\leq \sum_{\alpha \in \bQ^{\times}} \int_{T(\bA)}\left| \int_{N(\bA)} \sH (\mathbf s, g)^{-1} \psi(\alpha x) \, \mathrm dx \right| \, \mathrm da^{\times}\\
&= \sum_{\alpha \in \bQ^{\times}} \prod_p \int_{T(\bQ_p)}\left| \int_{N(\bQ_p)} \sH_p (\mathbf s, g_p)^{-1} \psi_p (\alpha x_p) \,\mathrm dx_p \right| \mathrm da_p^{\times}\\ 
&\qquad \qquad \qquad \times \int_{T(\bR)} \left| \int_{N(\bR)} \sH_\infty (\mathbf s, g_\infty)^{-1} \psi_\infty (\alpha x) \, \mathrm dx_\infty \right| \, \mathrm da_\infty^{\times}.
\end{align*}
Assume that $p \notin S$. Since the height function is right $\mathbf K_p$-invariant, we obtain that for any $y_p \in \bZ_p$,
\begin{align*}
\int_{N(\bQ_p)} \sH_p (\mathbf s, g_p)^{-1} \psi_p (\alpha x_p) \,\mathrm dx_p 
&= \int_{N(\bQ_p)} \sH_p (\mathbf s, (x_p + a_py_p, a_p))^{-1} \psi_p (\alpha x_p) \,\mathrm dx_p\\
&= \int_{N(\bQ_p)} \sH_p (\mathbf s, g_p)^{-1} \psi_p (\alpha x_p) \int_{\bZ_p} \overline{\psi}_p(\alpha a_p y_p)\, \mathrm dy_p \,\mathrm dx_p\\
&= 0 \qquad \text{ if $|\alpha a_p|_p >1$.}
\end{align*}
Hence we can conclude that
$$
\int_{T(\bQ_p)}\left| \int_{N(\bQ_p)} \sH_p (\mathbf s, g_p)^{-1} \psi_p (\alpha x_p) \,\mathrm dx_p \right| \mathrm da_p^{\times} \leq \int_{G(\bQ_p)} \sH_p (\Re (\mathbf s), g_p)^{-1} \boldsymbol{1}_{\bZ_p}(\alpha a_p) \, \mathrm dg_p.
$$
Similarly, for $p \in S$, we can conclude that
$$
\int_{T(\bQ_p)}\left| \int_{N(\bQ_p)} \sH_p (\mathbf s, g_p)^{-1} \psi_p (\alpha x_p) \,\mathrm dx_p \right| \mathrm da_p^{\times} \leq \int_{G(\bQ_p)} \sH_p (\Re (\mathbf s), g_p)^{-1} \boldsymbol{1}_{\frac{1}{N}\bZ_p}(\alpha a_p) \, \mathrm dg_p.
$$
Then the convergence of the following sum
$$
\sum_{\alpha \in \bQ^{\times}} \prod_p \int_{G(\bQ_p)} \sH_p^{-1} \boldsymbol{1}_{\frac{1}{N}\bZ_p}(\alpha a_p) \, \mathrm dg_p \cdot \int_{T(\bR)} \left| \int_{N(\bR)} \sH_\infty^{-1} \psi_\infty (\alpha x) \, \mathrm dx_\infty \right| \, \mathrm da_\infty^{\times},
$$
can be verified from the detailed study of the local integrals which we will conduct later. See proofs of Lemmas \ref{lemm:denominator}, \ref{lemm:real}, and \ref{lemm:Euler_product}.

Next we need to check that
$$
\int_{-\infty}^{+\infty} 
|(\sZ(\mathbf s, g), \theta_{m, \lambda, t})| \, \mathrm dt < +\infty.
$$
It is easy to see that
\begin{align*}
(\sZ(\mathbf s, g), \theta_{m, \lambda, t}) & = \int_{G(\bQ) \backslash G(\bA_{\bQ})} \sZ(\mathbf s, g) \overline{\theta}_{m, \lambda, t}\, \mathrm dg\\
&= \int_{G(\bA_{\bQ})} \sH (\mathbf s, g)^{-1} \overline{\theta}_{m, \lambda, t}\, \mathrm dg\\
&= \sum_{\alpha \in \bQ^{\times}} \int_{G(\bA_{\bQ})} \sH (\mathbf s, g)^{-1} \overline{\psi}(\alpha x)\overline{\mathbf v}_{m, \lambda}(\alpha a_{\text{fin}}) 
|\alpha a_{\infty}|_\infty^{-it}\, \mathrm dg\\
&= \sum_{\alpha \in \bQ^{\times}} \prod_p \sH_p'(\mathbf s, m, \lambda, \alpha) \cdot 
\sH_{\infty}'(\mathbf s, t, \alpha),
\end{align*}
where $\sH_p'(\mathbf s, m, \lambda, \alpha)$ 
is given by
\begin{align*}
&= \int_{G(\bQ_p)} \sH_p(\mathbf s, g_p)^{-1} 
\overline{\psi}_p(\alpha x_p) \boldsymbol{1}_{m \bZ_p^{\times}}(\alpha a_p)\, \mathrm dg_p, & p \notin S\\
&= \int_{G(\bQ_p)} \sH_p(\mathbf s, g_p)^{-1} \overline{\psi}_p(\alpha x_p) \overline{\lambda}_p(  \alpha a_p /p^{v_p(\alpha a_p)}) 
\boldsymbol{1}_{\frac{m}{N} \bZ_p^{\times}} (\alpha a_p)\, \mathrm dg_p, & p \in S
\end{align*}
and 
$$
\sH_\infty'(\mathbf s, t, \alpha)= 
\int_{G(\bR)} \sH_\infty(\mathbf s, g_\infty)^{-1} 
\overline{\psi}_\infty(\alpha x_\infty) |\alpha a_\infty|_\infty^{-it} \mathrm dg_\infty. 
$$
The integrability follows from the proof of Lemma \ref{lemm:real}. Thus we can apply Proposition \ref{prop:spectral}, and the identity in our statement follows from the continuity of $\sZ(\mathbf s, g)$.
\end{proof}

We obtained that
\begin{align*}
&\sZ_1(\mathbf s, \textnormal{id}) = \sum_{\lambda \in \mathbf M} \sum_{m=1}^{\infty} \frac{1}{4\pi} \int_{-\infty}^{+\infty} 
(\sZ(\mathbf s, g), \theta_{m, \lambda, t}) \theta_{m, \lambda, t}(\text{id})\, \mathrm dt\\
= &\sum_{\lambda \in \mathbf M, \, \lambda(-1) = 1} \sum_{m=1}^{\infty} \frac{1}{2\pi} \int_{-\infty}^{+\infty} 
(\sZ(\mathbf s, g), \theta_{m, \lambda, t}) 
\prod_{p \in S} \lambda_p \left(\frac{m}{N} \cdot p^{-v_p(m/N)} \right)\left| \frac{m}{N} \right|_\infty^{it} \,\mathrm dt.
\end{align*}

We will use the following notation:
\begin{align*}
\lambda_S(\alpha a_p) &:=  \prod_{q \in S} \overline{\lambda}_q(p^{v_p(\alpha a_p)}), &p \notin S\\
\lambda_{S, p}(\alpha a_p) &:= \overline{\lambda}_p\left(\frac{\alpha a_p}{p^{v_p(\alpha a_p)}}\right)\!\! \prod_{q \in S\setminus p}\!\!\! \overline{\lambda}_q(p^{v_p(\alpha a_p)}), &p \in S.
\end{align*}

\begin{prop}
\label{prop:spectral_general}
If $\Re (s)$ is sufficiently large, then
\begin{align*}
\sZ_1(\mathbf s, \text{id}) &= \sum_{\lambda \in \mathbf M, \, \lambda(-1) = 1} \sum_{\alpha \in \bQ^{\times}} \frac{1}{2\pi} \int_{-\infty}^{+\infty} \prod_p \widehat{\sH}_p (\mathbf s, \lambda, t, \alpha) \cdot \widehat{\sH}_{\infty}(\mathbf s, t, \alpha)\, \mathrm dt,
\end{align*}
where $\widehat{\sH}_p (\mathbf s, \lambda, t, \alpha)$ is given by 
\begin{align*}
&\int_{G(\bQ_p)}\!\!\!\! \!\!\!\!\sH_p(\mathbf s, g_p)^{-1} \overline{\psi}_p(\alpha x_p)\lambda_S(\alpha a_p) \boldsymbol{1}_{\bZ_p}(\alpha a_p) |a_p|_p^{-it} \mathrm dg_p,  &p \notin S\\
&\int_{G(\bQ_p)} \!\!\!\! \!\!\!\! \sH_p(\mathbf s, g_p)^{-1} \overline{\psi}_p(\alpha x_p) 
\lambda_{S, p}(\alpha a_p) \boldsymbol{1}_{\frac{1}{N} \bZ_p}(\alpha a_p) |a_p|_p^{-it} \mathrm dg_p,  & p \in S
\end{align*}
and 
$$
\widehat{\sH}_\infty(\mathbf s, t, \alpha)
= \int_{G(\bR)} \sH_\infty(\mathbf s, g_\infty)^{-1} 
\overline{\psi}_\infty(\alpha x_\infty) |a_\infty|_\infty^{-it}\, \mathrm dg_\infty
$$
\end{prop}
\begin{proof}
For simplicity, we assume that $S = \emptyset$. We have seen that
$$
\sZ_1(\mathbf s, \text{id})=\sum_{m=1}^{\infty} \sum_{\alpha \in \bQ^{\times}} \frac{1}{2\pi} \int_{-\infty}^{+\infty} \prod_p \sH_p'(\mathbf s, m, \alpha) \cdot 
\sH_{\infty}'(\mathbf s, t, \alpha) |m|_\infty^{it} \, \mathrm dt.
$$
On the other hand, it is easy to see that
$$
\widehat{\sH}_p (\mathbf s, t, \alpha) = \sum_{j = 0}^{\infty} \int_{G(\bQ_p)}\!\!\!\! \!\!\!\!\sH_p(\mathbf s, g_p)^{-1} \overline{\psi}_p(\alpha x_p) \boldsymbol{1}_{p^j\bZ_p^{\times}}(\alpha a_p) \left|\frac{p^j}{\alpha}\right|_p^{-it} \mathrm dg_p.
$$
Hence we have the formal identity:
$$
\prod_p \widehat{\sH}_p (\mathbf s, t, \alpha) \cdot \widehat{\sH}_\infty(\mathbf s, t, \alpha) = \sum_{m=1}^{\infty} \prod_p \sH_p'(\mathbf s, m, \alpha) \cdot 
\sH_{\infty}'(\mathbf s, t, \alpha) |m|_\infty^{it},
$$
and our assertion follows from this. To justify the above identity, 
we need to address convergence issues; this will be discussed below 
(see the proof of Lemma~\ref{lemm:denominator}).
\end{proof}

Thus we need to study the local integrals in Proposition \ref{prop:spectral_general}. We introduce some notation:
\begin{align*}
&\mathcal I_1 = \{ \iota \in \mathcal I \, |\, D_\iota \subset \text{Supp}(\text{div}_0(a)) \}\\
&\mathcal I_2 = \{ \iota \in \mathcal I \, |\, D_\iota \subset \text{Supp}(\text{div}_\infty(a)) \}\\
&\mathcal I_3 = \{ \iota \in \mathcal I \, |\, D_\iota \not \subset \text{Supp}(\text{div}(a)) \}.
\end{align*}
Note that $\mathcal I = \mathcal I_1 \sqcup \mathcal I_2 \sqcup \mathcal I_3$ and 
$\mathcal I_1 \neq \emptyset$. Also $D_\iota \subset \text{Supp}(\text{div}_\infty(x))$ for any $\iota \in \mathcal I_3$ because 
$$D = \cup_{\iota \in \mathcal I} D_\iota = 
\text{Supp}(\text{div}(a)) \cup \text{Supp}(\text{div}_\infty(x)).
$$ 
Let
$$
-\text{div}(\omega) = \sum_{\iota \in \mathcal I} d_\iota D_\iota,
$$
where $\omega = \mathrm dx \mathrm da /a$ is the top degree right invariant form on $G$. Note that $\omega$ defines a measure $|\omega|$ on an analytic manifold $G(\bQ_v)$, and for any finite place $p$, 
$$
|\omega| = \left(1-\frac{1}{p}\right)\mathrm dg_p,
$$
where $\mathrm dg_p$ is the standard Haar measure defined in Section \ref{sect:harm}.
 
\begin{lemm}
\label{lemm:domain}
Consider an open convex cone $\Omega$ in $\Pic^G (X)_\bR$, defined by the following relations:
$$
\begin{cases}
s_\iota - d_\iota +1> 0 & \textnormal{ if } \iota \in \mathcal I_1\\
s_\iota -d_\iota +1 + e_\iota> 0 & \textnormal{ if } \iota \in \mathcal I_2\\
s_\iota -d_\iota +1 > 0 & \textnormal{ if } \iota \in \mathcal I_3\\
\end{cases}
$$
where $e_\iota = |\textnormal{ord}_{D_\iota}(x)|$. 
Then $\widehat{\sH}_p (\mathbf s, \lambda, t, \alpha)$ and 
$\widehat{\sH}_\infty(\mathbf s, t, \alpha)$  are holomorphic on $\mathsf T_\Omega$.
\end{lemm}
\begin{proof}
First we prove our assertion for $\widehat{\sH}_\infty$. 
We can assume that 
$$
\widehat{\sH}_v(\mathbf s, t) = \widehat{\sH}_v(\mathbf s - it\mathbf m(a), 0),
$$
where $\mathbf m(a)\in \mathfrak X^*(G)\subset \Pic^G(X)$ is the character associated to the rational function $a$
(by choosing an appropriate height function).
It suffices to discuss the case when $t=0$.
Choose a finite covering $\{ U_\eta \}$ of $X(\bR)$ by open subsets
and local coordinates $y_\eta, z_\eta$ on $U_\eta$ 
such that the union of the boundary divisor $D$ and $E$ is 
locally defined by $y_\eta =0$ or $y_\eta \cdot z_\eta = 0$.
Choose a partition of unity $\{ \theta_\eta \}$;
the local integral takes the form
$$
\widehat{\sH}_\infty(\mathbf s, \alpha) = \sum_{\eta} \int_{G(\bR)} \sH_\infty(\mathbf s, g_\infty)^{-1} 
\overline{\psi}_\infty(\alpha x_\infty) \theta_\eta \, \mathrm dg_\infty. 
$$
Each integral is a oscillatory integral in the variables $y_\eta, z_\eta$.
For example, assume that $U_\eta$ meets $D_{\iota}, D_{\iota'}$, 
where $\iota, \iota' \in \mathcal I_2$.
Then
\begin{multline*}
\int_{G(\bR)} \sH_\infty(\mathbf s, g_\infty)^{-1} 
\overline{\psi}_\infty(\alpha x_\infty) \theta_\eta \mathrm dg_\infty \\
= \int_{\bR^2} |y_\eta|^{s_{\iota} - d_{\iota}} |z_\eta|^{s_{\iota'} - d_{\iota'}} 
\overline{\psi}_\infty\left(\frac{\alpha f}{y_\eta^{e_{\iota}}z_\eta^{e_{\iota'}}}\right) \Phi( \mathbf s, y_\eta, z_\eta) \mathrm d y_\eta \,\mathrm d z_\eta,
\end{multline*}
where $\Phi$ is a smooth function with compact support 
and $f$ is a nonvanishing analytic function.
Shrinking $U_\eta$ and changing variables, if necessary, we may assume that $f$ is a constant.
Proposition \ref{prop:osillatory} implies that this integral is holomorphic everywhere.
The other integrals can be studied similarly.

Next we consider finite places. 
Let $p$ be a prime of good reduction. 
Since
$$
\text{Supp}(\text{div}_0(a)) \cap \text{Supp}(\text{div}_\infty(a)) =  \emptyset,
$$
the smooth function $\boldsymbol 1_{\bZ_p} (\alpha a_p)$ extends 
to a smooth function $h$ on $X(\bQ_p)$.
Let 
$$
U =\{ h = 1\}.
$$
Then
$$
\widehat{\sH}_p (\mathbf s, \lambda, \alpha) = \int_{U} \sH_p(\mathbf s, g_p)^{-1}
 \overline{\psi}_p(\alpha x_p) \lambda_S(\alpha a_p) \mathrm dg_p.
$$
Now the proof of \cite[ Lemma 4.4.1]{volume} 
implies that this is holomorphic on $\mathsf T_\Omega$ 
because $U \cap (\cup_{\iota \in \mathcal I_2} D_\iota(\bQ_p)) = \emptyset$.
Places of bad reduction are treated similarly.
\end{proof}

\begin{lemm}
\label{lemm:denominator}
Let $|\alpha|_p = p^k >1$. Then, for any compact set in $\Omega$ and for any $\delta > 0$, there exists a constant $C>0$ such that
$$
|\widehat{\sH}_p (\mathbf s, \lambda, t, \alpha)| < C |\alpha|_p^{-\min_{\iota \in \mathcal I_1} \{\Re (s_\iota) - d_\iota +1- \delta \}},
$$
for $\Re(\mathbf s)$ in that compact set.
\end{lemm}
\begin{proof}
First assume that $p$ is a good reduction place. 
Let $\rho : \mathcal X(\bZ_p) \ra \mathcal X (\bF_p)$ be the reduction map modulo $p$
where $\mathcal X$ is a smooth integral model of $X$ over $\Spec(\bZ_p)$.
Note that 
$$
\rho(\{ |a|_p <1 \}) \subset \cup_{\iota \in \mathcal I_1} 
\mathcal D_\iota(\bF_p),
$$ 
where $\mathcal D_\iota$ is the Zariski closure of $D_\iota$ in $\mathcal X$.
Thus $\widehat{\sH}_p (\mathbf s, \lambda, \alpha)$ is given by
$$
\widehat{\sH}_p (\mathbf s, \lambda, \alpha)
= \sum_{\tilde{x} \in \cup_{\iota \in \mathcal I_1} \mathcal D_\iota(\bF_p)} \int_{\rho^{-1}(\tilde{x})} \sH_p(\mathbf s, g_p)^{-1}  \overline{\psi}_p(\alpha x_p) \lambda_S(\alpha a_p) \boldsymbol{1}_{\bZ_p}(\alpha a_p) \mathrm dg_p.
$$
Let $\tilde{x} \in \mathcal D_\iota(\bF_p)$ for some $\iota \in \mathcal I_1$, but $\tilde{x} \notin \mathcal D_{\iota'}(\bF_p)$ for any $\iota' \in \mathcal I \setminus \{ \iota \}$. Since $p$ is a good reduction place, we can find analytic coordinates $y, z$ such that
\begin{align*}
\left| \int_{\rho^{-1}(\tilde{x})} \right| 
&\leq \int_{\rho^{-1}(\tilde{x})} \sH_p(\Re (\mathbf s), g_p)^{-1} \boldsymbol{1}_{\bZ_p}(\alpha a_p) \mathrm dg_p\\
&= \left( 1 - \frac{1}{p} \right)^{-1}\int_{\rho^{-1}(\tilde{x})} \sH_p(\Re (\mathbf s) - \mathbf d, g_p)^{-1} \boldsymbol{1}_{\bZ_p}(\alpha a_p) \mathrm d\tau_{X,p}\\
&= \left( 1 - \frac{1}{p} \right)^{-1}\int_{\mm_p^2} |y|_p^{\Re (s_\iota) - d_\iota} \boldsymbol{1}_{\bZ_p}(\alpha y) \mathrm dy_p \mathrm dz_p\\
&= \frac{1}{p} \cdot \frac{p^{-k(\Re (s_\iota) - d_\iota +1)}}{1 - p^{-(\Re (s_\iota) - d_\iota +1)}},
\end{align*}
where $\mathrm d\tau_{X,p}$ is the local Tamagawa measure 
(see \cite[Section 2]{volume} for the definition).
For the construction of such local analytic coordinates, 
see \cite{Weil}, \cite{denef}, or \cite{salberger}.
If $\tilde{x} \in \mathcal D_\iota(\bF_p) \cap \mathcal D_{\iota'}(\bF_p)$ for $\iota \in \mathcal I_1$, $\iota' \in \mathcal I_3$, then we can find local 
analytic coordinates $y, z$ such that
\begin{align*}
\left| \int_{\rho^{-1}(\tilde{x})} \right|
&\leq \left( 1 - \frac{1}{p} \right) \int_{\mm_p^2} |y|_p^{\Re (s_\iota) - d_\iota + 1} |z|_p^{\Re (s_{\iota'}) - d_{\iota'} + 1}  \boldsymbol{1}_{\bZ_p}(\alpha y) \mathrm dy^{\times}_p \mathrm dz^{\times}_p\\
&= \left( 1 - \frac{1}{p} \right) \frac{p^{-k(\Re (s_\iota) - d_\iota +1)}}{1 - p^{-(\Re (s_\iota) - d_\iota +1)}} \frac{p^{-(\Re (s_{\iota'}) - d_{\iota'} +1)}}{1 - p^{-(\Re (s_{\iota'}) - d_{\iota'} +1)}}.
\end{align*}
If $\tilde{x} \in \mathcal D_\iota(\bF_p) \cap \mathcal D_{\iota'}(\bF_p)$ for $\iota, \iota' \in \mathcal I_1$, $\iota \neq \iota'$, then we can find analytic coordinates $x, y$ such that
\begin{align*}
\left| \int_{\rho^{-1}(\tilde{x})} \right|
&\leq \left( 1 - \frac{1}{p} \right) \int_{\mm_p^2} |y|_p^{\Re (s_\iota) - d_\iota + 1} |z|_p^{\Re (s_{\iota'}) - d_{\iota'} + 1}  \boldsymbol{1}_{\bZ_p}(\alpha yz) \, \mathrm dy^{\times}_p \mathrm dz^{\times}_p\\
&\leq \left( 1 - \frac{1}{p} \right) \int_{\mm_p^2} |yz|_p^{\min \{\Re (s_\iota) - d_\iota + 1, \, \Re (s_{\iota'}) - d_{\iota'} + 1 \}}  \boldsymbol{1}_{\bZ_p}(\alpha yz) \, \mathrm dy^{\times}_p \mathrm dz^{\times}_p\\
&= \left( 1 - \frac{1}{p} \right) \left( (k-1) \frac{p^{-kr}}{1-p^{-r}} + \frac{p^{-(k+1)r}}{(1-p^{-r})^2} \right),
\end{align*}
where 
$$
r = \min \{\Re (s_\iota) - d_\iota + 1, \, \Re (s_{\iota'}) - d_{\iota'} + 1 \}.
$$ 
It follows from these inequalities and Lemma 9.4 in \cite{VecIII} 
that there exists a constant $C> 0$, independent of $p$, 
satisfying the inequality in the statement.

Next assume that $p$ is a bad reduction place. Choose 
an open covering $\{U_\eta\}$ of $\cup_{\iota \in \mathcal I_1} D_\iota(\bQ_p)$ such that 
$$
(\cup_\eta U_\eta) \cap (\cup_{\iota \in \mathcal I_2} D_\iota(\bQ_p))
 = \emptyset,
$$ 
and each $U_\eta$ has analytic coordinates $y_\eta, z_\eta$. Moreover, we can assume that the boundary divisor is defined by $y_\eta = 0$ or $y_\eta \cdot z_\eta = 0$ on $U_\eta$. Let $V$ be the complement of $\cup_{\iota \in \mathcal I_1} D_\iota(\bQ_p)$, and consider the partition of unity for $\{ U_\eta, V \}$ which we denote by $\{ \theta_\eta, \theta_V \}$. If $k$ is sufficiently large, then
$$
\{ \boldsymbol{1}_{\frac{1}{N}\bZ_p}(\alpha a) = 1 \} \cap \text{Supp}(\theta_V) = \emptyset .
$$
Hence if $k$ is sufficiently large, then 
\begin{align*}
|\widehat{\sH}_p (\mathbf s, \lambda, \alpha)|
&\leq \sum_\eta \int_{U_\eta}  \sH_p(\Re (\mathbf s), g_p)^{-1} \boldsymbol{1}_{\frac{1}{N}\bZ_p}(\alpha a_p)\cdot \theta_\eta \, \mathrm dg_p.
\end{align*}
When $U_\eta$ meets 
only one component $D_\iota(\bQ_p)$ for $\iota \in \mathcal I_1$, then
\begin{align*}
\int_{U_\eta} &\leq \int_{\bQ_p^2} |y_\eta|_p^{\Re (s_\iota) - d_\iota} \boldsymbol{1}_{c\bZ_p}(\alpha y_\eta) \Phi ( \mathbf s, y_\eta, z_\eta) \,\mathrm d y_{\eta,p} 
\mathrm d z_{\eta,p} \ll p^{-k(\Re (s_\iota) - d_\iota +1)},
\end{align*}
as $k \ra \infty$, where $c$ is some rational number and $\Phi$ 
is a smooth function with compact support.
Other integrals are treated similarly. 
\end{proof}

We record the following useful lemma (see, e.g., \cite[Lemma 2.3.1]{CT-additive}):
\begin{lemm}
\label{lemm:usefulII}
Let $d$ be a positive integer and $a \in \bQ_p$. If $|a|_p >p$ and 
$p\nmid d$, then
$$
\int_{\bZ_p^{\times}} \overline{\psi}_p(a x^d) \,\mathrm dx^{\times}_p = 0.
$$
Moreover, if $|a|_p = p$ and $d = 2$, then
$$
\int_{\bZ_p^{\times}} \overline{\psi}_p(a x^d) \mathrm dx^{\times}_p= 
\begin{cases}
\frac{\sqrt{p}-1}{p-1} \text{ or } \frac{i\sqrt{p} - 1}{p-1} & \textnormal{ if 
$pa$ is a quadratic residue} ,\\
\frac{-\sqrt{p}-1}{p-1} \text{ or } \frac{-i\sqrt{p} - 1}{p-1} & \textnormal{ if $pa$ is a quadratic non-residue.} 
\end{cases}
$$
\end{lemm}

\begin{lemm}
\label{lemm:numerator}
Let $|\alpha|_p =p^{-k}<1$. 
Consider an open convex cone $\Omega_\epsilon$ in $\Pic (X)_\bR$, defined by the following relations:
$$
\begin{cases}
s_\iota - d_\iota +1> 0 & \textnormal{ if } \iota \in \mathcal I_1\\
s_\iota -d_\iota +2+\epsilon> 0 & \textnormal{ if } \iota \in \mathcal I_2\\
s_\iota -d_\iota +1 > 0 & \textnormal{ if } \iota \in \mathcal I_3\\
\end{cases}
$$
where $0 < \epsilon < 1/3$.
Then, for any compact set in $\Omega_\epsilon$, there exists a constant $C>0$ such that
$$
|\widehat{\sH}_p (\mathbf s, \lambda, t, \alpha)| < C
|\alpha|_p^{-\frac{2}{3}(1+2\epsilon)},
$$
for $\Re(\mathbf s)$ in that compact set.
\end{lemm}
\begin{proof}
First assume that $p$ is a good reduction place and that 
$p\nmid e_\iota$, for any $\iota \in \mathcal I_2$. We have
$$
\widehat{\sH}_p (\mathbf s, \lambda, \alpha)
= \sum_{\tilde{x} \in \mathcal X(\bF_p)} \int_{\rho^{-1}(\tilde{x})} \sH_p(\mathbf s, g_p)^{-1}\overline{\psi}_p(\alpha x_p) \lambda_S(\alpha a_p) \boldsymbol{1}_{\bZ_p}(\alpha a_p)  \, \mathrm dg_p.
$$
A formula of J. Denef 
(see \cite[Theorem 3.1]{denef} or \cite[Proposition 4.1.7]{volume}) 
and Lemma 9.4 in \cite{VecIII} give 
us an uniform bound:
$$
|\sum_{\tilde{x} \notin \cup_{\iota \in \mathcal I_2} \mathcal D_\iota(\bF_p)}|
\leq \sum_{\tilde{x} \notin \cup_{\iota \in \mathcal I_2}\mathcal D_\iota(\bF_p)} \int_{\rho^{-1}(\tilde{x})} \sH_p(\Re (\mathbf s), g_p)^{-1}\,\mathrm dg_p.
$$
Hence we need to study 
$$
\sum_{\tilde{x} \in \cup_{\iota \in \mathcal I_2} \mathcal D_\iota(\bF_p)} \int_{\rho^{-1}(\tilde{x})} \sH_p(\mathbf s, g_p)^{-1} 
\overline{\psi}_p(\alpha x_p) \lambda_S(\alpha a_p) \boldsymbol{1}_{\bZ_p}(\alpha a_p)  \,\mathrm dg_p.
$$
Let $\tilde{x} \in \mathcal D_\iota (\bF_p)$ for some $\iota \in \mathcal I_2$, but $\tilde{x} \notin \mathcal D_{\iota'}(\bF_p) \cup \mathcal E(\bF_p)$ for any $\iota' \in \mathcal I \setminus \{ \iota \}$, where $\mathcal E$ is the Zariski closure of $E$ in $\mathcal X$. Then we can find local analytic coordinates $y, z$ such that
$$
\int_{\rho^{-1}(\tilde{x})} = \left( 1 - \frac{1}{p} \right)^{-1} \int_{\mm_p^2} |y|_p^{s_\iota - d_\iota}  \overline{\psi}_p(\alpha f/y^{e_\iota}) \lambda_S(\alpha y^{-1}) \boldsymbol{1}_{\bZ_p}(\alpha y^{-1}) \,\mathrm dy_p \mathrm dz_p,
$$
where $f \in \bZ_p [[y, z]]$ such that $f(0) \in \bZ_p^{\times}$. Since $p$ does not divide $e_\iota$, there exists $g \in \bZ_p[[y, z]]$ such that $f = f(0)g^{e_\iota}$. After a change of variables, 
we can assume that $f = u \in \bZ_p^{\times}$. Lemma \ref{lemm:usefulII} implies that
\begin{align*}
\int_{\rho^{-1}(\tilde{x})} &= \frac{1}{p} \int_{\mm_p} |y|_p^{s_\iota - d_\iota + 1} \lambda_S(\alpha y^{-1}) \int_{\bZ_p^{\times}} \overline{\psi}_p(\alpha ub^{e_\iota}/y^{e_\iota}) \mathrm db^{\times}_p \boldsymbol{1}_{\bZ_p}(\alpha y^{-1}) 
\mathrm dy^{\times}_p\\
&= \frac{1}{p} \int_{p^{-(k+1)} \leq |y^{e_\iota}|_p} 
|y|_p^{s_\iota - d_\iota + 1} \lambda_S(\alpha y^{-1}) 
\int_{\bZ_p^{\times}} \overline{\psi}_p(\alpha ub^{e_\iota}/y^{e_\iota})\, 
\mathrm db^{\times}_p \mathrm dy^{\times}_p
\end{align*}
Thus it follows from the second assertion of Lemma \ref{lemm:usefulII} that
\begin{align*}
\left| \int_{\rho^{-1}(\tilde{x})} \right| &\leq 
\frac{1}{p} \int_{p^{-(k+1)} \leq |y^{e_\iota}|} |y|_p^{\Re (s_\iota) - d_\iota + 1} \left| \int_{\bZ_p^{\times}} \overline{\psi}_p(\alpha ub^{e_\iota}/y^{e_\iota}) \mathrm db^{\times}_p \right| \mathrm dy^{\times}_p\\
&\leq\frac{1}{p}k p^{\frac{k}{e_\iota}(1+\epsilon)} + \frac{1}{p}p^{\frac{k+1}{e_\iota}(1+\epsilon)} \times
\begin{cases}
1 & \text{ if } e_\iota >2\\
\frac{1}{\sqrt{p} -1} & \text{ if } e_\iota =2
\end{cases}\\
&\ll \frac{1}{p}k p^{\frac{2}{3}k(1+\epsilon)}.
\end{align*}
If $\tilde{x} \in \mathcal D_\iota (\bF_p) \cap \mathcal E(\bF_p)$, for some $\iota \in \mathcal I_2$, then we have
\begin{align*}
\int_{\rho^{-1}(\tilde{x})} &= \left( 1 - \frac{1}{p} \right)^{-1} 
\int_{\mm_p^2} |y|_p^{s_\iota - d_\iota} \overline{\psi}_p(\alpha z/y^{e_\iota}) \lambda_S(\alpha y^{-1}) \boldsymbol{1}_{\bZ_p}(\alpha y^{-1}) \, 
\mathrm dy_p \mathrm dz_p\\
&= \int_{\mm_p} |y|_p^{s_\iota - d_\iota + 1 } \lambda_S(\alpha y^{-1}) \boldsymbol{1}_{\bZ_p}(\alpha y^{-1}) \int_{\mm_p} \overline{\psi}_p(\alpha z/y^{e_\iota}) \mathrm dz_p \mathrm dy^{\times}_p\\
&= \frac{1}{p} \int_{p^{-(k+1)} \leq |y|_p^{e_\iota} <1} |y|_p^{s_\iota - d_\iota + 1 } \lambda_S(\alpha y^{-1}) \, \mathrm dy^{\times}_p.
\end{align*}
Hence we obtain that
$$
\left| \int_{\rho^{-1}(\tilde{x})} \right| \leq \frac{1}{p} \int_{p^{-(k+1)} 
\leq |y|_p^{e_\iota} <1} |y|_p^{\Re (s_\iota) - d_\iota + 1 } \,\mathrm dy^{\times}_p 
\leq kp^{\frac{k}{e_\iota}(1+\epsilon)} < kp^{\frac{2}{3}k(1+\epsilon)}.
$$
If $\tilde{x} \in \mathcal D_\iota (\bF_p) \cap \mathcal D_{\iota'} (\bF_p)$ for some $\iota \in \mathcal I_2$ and $\iota' \in \mathcal I_3$, then it follows from Lemma \ref{lemm:usefulII}
\begin{align*}
\int_{\rho^{-1}(\tilde{x})}&= \left( 1 - \frac{1}{p} \right)^{-1} 
\int_{\mm_p^2} |y|_p^{s_\iota - d_\iota} |z|_p^{s_{\iota'} - d_{\iota'}} 
\overline{\psi}_p \left(\frac{\alpha u}{y^{e_\iota}z^{e_{\iota'}}}\right) \lambda_S(\alpha y^{-1}) \boldsymbol{1}_{\bZ_p}(\alpha y^{-1}) 
\, \mathrm dy_p \mathrm dz_p\\
&= \left( 1 - \frac{1}{p} \right)^{-1} 
\int 
|y|_p^{s_\iota - d_\iota} |z|_p^{s_{\iota'} - d_{\iota'}} 
\lambda_S(\alpha y^{-1}) \int_{\bZ_p^{\times}} \overline{\psi}_p \left(\frac{\alpha ub^{e_\iota}}{y^{e_\iota}z^{e_{\iota'}}}\right)\, 
\mathrm db^{\times}_p \, \mathrm dy_p \mathrm dz_p,
\end{align*}
where the last integral is over the domain
$$
\{ (y, z) \in \mm_p^2 : p^{-(k+1)} \leq  |y^{e_{\iota}} z^{e_{\iota'}}|_p \}.
$$
We conclude that
\begin{align*}
\left| \int_{\rho^{-1}(\tilde{x})} \right| &\leq
\left( 1 - \frac{1}{p} \right)^{-1} 
\int_{p^{-(k+1)} \leq |y^{e_\iota}z^{e_{\iota'}}|_p} 
|y|_p^{\Re(s_\iota) - d_\iota} |z|_p^{\Re(s_{\iota'}) - d_{\iota'}}\, \mathrm dy_p 
\mathrm dz_p
\\
&\leq \left( 1 - \frac{1}{p} \right)^{-1} 
\int_{p^{-k} \leq |y^{e_\iota}|_p <1} |y|_p^{\Re(s_\iota) - d_\iota} 
\, \mathrm dy_p \int_{\mm_p} |z|_p^{\Re(s_{\iota'}) - d_{\iota'}}\, \mathrm dz_p\\
&\leq k p^{\frac{k}{e_\iota}(1+\epsilon)}
\frac{p^{-(\Re(s_{\iota'}) - d_{\iota'} +1)}}{1 - p^{-(\Re(s_{\iota'}) - d_{\iota'} +1)}}.
\end{align*}
If $\tilde{x} \in \mathcal D_\iota (\bF_p) \cap \mathcal D_{\iota'} (\bF_p)$ for some $\iota, \iota' \in \mathcal I_2$, then the local integral on $\rho^{-1}(\tilde{x})$ is:
\begin{align*}
\left( 1 - \frac{1}{p} \right)^{-1} 
\int_{\mm_p^2} |y|_p^{s_\iota - d_\iota} |z|_p^{s_{\iota'} - d_{\iota'}} 
\overline{\psi}_p \left(\frac{\alpha u}{y^{e_\iota}z^{e_{\iota'}}}\right) \lambda_S(\alpha y^{-1}z^{-1})\boldsymbol{1}_{\bZ_p}(\alpha y^{-1}z^{-1}) 
 \,\mathrm dy_p \mathrm dz_p\\
= \left( 1 - \frac{1}{p} \right) \int_{\mm_p^2} 
|y|_p^{s_\iota - d_\iota} |z|_p^{s_{\iota'} - d_{\iota'}} 
\lambda_S(\alpha y^{-1}z^{-1}) \int_{\bZ_p^{\times}} 
\overline{\psi}_p \left(\frac{\alpha ub^{e_\iota}}{y^{e_\iota}z^{e_{\iota'}}}\right)
\, \mathrm db^{\times}_p \mathrm dy^{\times}_p \mathrm dz^{\times}_p.
\end{align*}
We can assume that $e_\iota \leq e_{\iota'}$. 
Then we can conclude that
\begin{align*}
\left| \int_{\rho^{-1}(\tilde{x})} \right| &\leq 
\int_{p^{-k} \leq |y^{e_\iota}z^{e_{\iota'}}|_p} 
|y^{e_\iota}z^{e_{\iota'}}|_p^{-\frac{1}{e_\iota}(1+\epsilon)}  \,
\mathrm dy^{\times}_p \mathrm dz^{\times}_p\\ 
&+ \int_{p^{-(k+1)} = |y^{e_\iota}z^{e_{\iota'}}|_p} 
|y^{e_\iota}z^{e_{\iota'}}|_p^{-\frac{1}{e_\iota}(1+\epsilon)} \left| \int_{\bZ_p^{\times}} \overline{\psi}_p \left(\frac{\alpha ub^{e_\iota}}{y^{e_\iota}z^{e_{\iota'}}}\right) \mathrm db^{\times} \right| \,\mathrm dy^{\times}_p \mathrm dz^{\times}_p\\
&\leq k^2 p^{\frac{k}{e_\iota}(1+\epsilon)} + kp^{\frac{k+1}{e_\iota}(1+\epsilon)} \times
\begin{cases}
1 & \text{ if } e_\iota >2\\
\frac{1}{\sqrt{p} -1} & \text{ if } e_\iota =2
\end{cases}\\
&\ll k^2 p^{\frac{2}{3}k(1+\epsilon)}.
\end{align*}
Thus our assertion follows from these estimates and Lemma 9.4 in \cite{VecIII}.

Next assume that $p$ is a place of bad reduction or that $p$ divides 
$e_\iota$, for some $\iota \in \mathcal I_2$. Fix a compact subset of $\Omega_\epsilon$ and assume that $\Re (\mathbf s)$ is in that compact set. Choose a finite open covering $\{ U_\eta \}$ of $\cup_{\iota \in \mathcal I_2} D_\iota(\bQ_p)$ with analytic coordinates $y_\eta, z_\eta$ such that the union of the boundary $D(\bQ_p)$ and $E(\bQ_p)$ is defined by $y_\eta = 0$ or $y_\eta \cdot z_\eta =0$. Let $V$ be the complement of $\cup_{\iota \in \mathcal I_2} D_\iota(\bQ_p)$, and consider a partition of unity $\{ \theta_\eta, \theta_V \}$ for $\{U_\eta, V\}$. Then it is clear that
$$
\int_{V} \sH_p(\mathbf s, g_p)^{-1} \overline{\psi}_p(\alpha x_p) 
\lambda_{S, p}(\alpha a_p) \boldsymbol{1}_{\frac{1}{N} \bZ_p}(\alpha a_p) \theta_V \mathrm dg_p,
$$
is bounded, so we need to study
$$
\int_{U_\eta} \sH_p(\mathbf s, g_p)^{-1} \overline{\psi}_p(\alpha x_p) 
\lambda_{S, p}(\alpha a_p) \boldsymbol{1}_{\frac{1}{N} \bZ_p}(\alpha a_p) \theta_{U_\eta} \mathrm dg_p.
$$
Assume that $U_\eta$ meets only one $D_\iota(\bQ_p)$ for some $\iota \in \mathcal I_2$. Then, the above integral looks like
$$
\int_{U_\eta} = \int_{\bQ_p^2} |y_\eta|_p^{s_\iota -d_\iota} \overline{\psi}_p(\alpha f/y_\eta^{e_\iota}) ) \lambda_{S, p}(\alpha g/y_\eta) \boldsymbol{1}_{\frac{1}{N} \bZ_p}(\alpha g/y_\eta) \Phi(\mathbf s, y_\eta, z_\eta) \mathrm dy_{\eta, p} \mathrm dz_{\eta, p},
$$
where $f$ and $g$ are nonvanishing analytic functions, and $\Phi$ is a smooth function with compact support. By shrinking $U_\eta$ and changing variables, if necessary, we can assume that $f$ and $g$ are constant. The proof of Proposition \ref{prop:osillatorybounds} implies our assertion for this integral. Other integrals are treated similarly.
\end{proof}

\begin{lemm}
\label{lemm:real}
For any compact set in an open convex cone $\Omega'$, defined by
$$
\begin{cases}
s_\iota - d_\iota -1 > 0  & \textnormal{ if } \iota \in \mathcal I_1\\
s_\iota - d_\iota +3> 0 & \textnormal{ if } \iota \in \mathcal I_2\\
s_\iota -d_\iota +1> 0 & \textnormal{ if } \iota \in \mathcal I_3\\
\end{cases}
$$
there exists a constant $C>0$ such that
$$
|\widehat{\sH}_\infty(\mathbf s, t, \alpha)| < \frac{C}{|\alpha|^2(1+t^2)},
$$
for $\Re(\mathbf s)$ in that compact set.
\end{lemm}
\begin{proof}
Consider the left invariant differential operators $\partial_a = a\partial /\partial a$ and $\partial_x = a\partial/\partial x$. Assume that $\Re (\mathbf s) \gg 0$. Integrating by parts, we have
\begin{align*}
\widehat{\sH}_\infty(\mathbf s, t, \alpha)
&= -\frac{1}{t^2} \int_{G(\bR)} \partial_a^2 \sH_\infty(\mathbf s, g_\infty)^{-1} 
\overline{\psi}_\infty(\alpha x_\infty) |a_\infty|_\infty^{-it}\, \mathrm dg_\infty\\
&= \frac{1}{(2\pi)^2 |\alpha|^2t^2} \int_{G(\bR)} \frac{\partial^2}{\partial x^2} (\partial_a^2 \sH_\infty(\mathbf s, g_\infty)^{-1}) 
\overline{\psi}_\infty(\alpha x_\infty) |a_\infty|_\infty^{-it}\, \mathrm dg_\infty.
\end{align*}
According to Proposition 2.2. in \cite{VecIII}, 
\begin{align*}
\frac{\partial^2}{\partial x^2} (\partial_a^2 \sH_\infty(\mathbf s, g_\infty)^{-1})
&= |a|^{-2} \partial_x^2 \partial_a^2 \sH_\infty(\mathbf s, g_\infty)^{-1}\\
&= \sH_\infty(\mathbf s - 2\mathbf m(a), g_\infty)^{-1} \times ( \text{a bounded smooth function}).
\end{align*}
Moreover, Lemma 4.4.1. of \cite{volume} tells us that
$$
\int_{G(\bR)} \sH_\infty(\mathbf s - 2\mathbf m(a), g_\infty)^{-1} \mathrm dg_\infty,
$$
is holomorphic on $\mathsf T_{\Omega'}$. Thus we can conclude our lemma.
\end{proof}

\begin{lemm}
\label{lemm:Euler_product}
The Euler product
$$
\prod_p \widehat{\sH}_p(\mathbf{s}, \lambda, t, \alpha) \cdot \widehat{\sH}_\infty (\mathbf{s}, t, \alpha)
$$
is holomorphic on $\mathsf T_{\Omega'}$.
\end{lemm}
\begin{proof}
First we prove that the Euler product is holomorphic on $\mathsf T_{\Omega'}$. To conclude this, we only need to discuss:
$$
\prod_{p \notin S\cup S_3, \, |\alpha|_p=1, }  \widehat{\sH}_p(\mathbf{s}, \lambda, t, \alpha),
$$
where $S_3 = \{ p : p\mid e_\iota \text{ for some } \iota \in \mathcal I_3\}$. Let $p$ be a prime such that $p \notin S\cup S_3$ and $|\alpha|_p =1$. Fix a compact subset of $\Omega'$, and assume that $\Re(s)$ is sitting in that compact set. From the definition of $\Omega'$, there exists $\epsilon > 0$ such that
$$
\begin{cases}
s_\iota - d_\iota + 1 > 2+\epsilon & \text{ for any } \iota \in \mathcal I_1\\
s_\iota - d_\iota + 1 > \epsilon & \text{ for any } \iota \in \mathcal I_3.
\end{cases}
$$
Since we have
$$
\{ |a|_p \leq 1 \} = X(\bQ_p) \setminus \rho^{-1}(\cup_{\iota \in \mathcal I_2} \mathcal D_\iota(\bF_p)),
$$
we can conclude that
$$
\widehat{\sH}_p(\mathbf{s}, \lambda, \alpha)
= \sum_{\tilde{x} \notin \cup_{\iota \in \mathcal I_2} \mathcal D_\iota(\bF_p)}
\int_{\rho^{-1}(\tilde{x})} \sH_p(\mathbf s, g_p)^{-1} \overline{\psi}_p(\alpha x_p) \lambda_S(a_p) \mathrm dg_p.
$$
It is easy to see that
$$
\sum_{\tilde{x} \notin \cup_{\iota \in \mathcal I} \mathcal D_\iota(\bF_p) }\int_{\rho^{-1}(\tilde{x})} 
= \int_{G(\bZ_p)} 1\, \mathrm dg_p =1.
$$
Also it follows from a formula of J. Denef (see \cite[Theorem 3.1]{denef} or \cite[Proposition 4.1.7]{volume}) and Lemma 9.4 in \cite{VecIII} that there exists an uniform bound $C > 0$ such that for any 
$
\tilde{x} \in \cup_{\iota \in \mathcal I_1} \mathcal D_\iota(\bF_p),
$
\begin{align*}
\left|\int_{\rho^{-1}(\tilde{x})} \right| 
< \int_{\rho^{-1}(\tilde{x})} \sH_p(\Re (\mathbf s), g_p)^{-1} \mathrm dg_p < \frac{C}{p^{2+\epsilon}}.
\end{align*}
Hence we need to obtain uniform bounds of $\int_{\rho^{-1}(\tilde{x})}$ for 
$$
\tilde{x} \in \cup_{\iota \in \mathcal I_3} \mathcal D_\iota(\bF_p) \setminus \cup_{\iota \in \mathcal I_1 \cup \mathcal I_2} \mathcal D_\iota(\bF_p).
$$
Let $\tilde{x} \in \mathcal D_\iota(\bF_p)$ for some $\iota \in \mathcal I_3$, but $\tilde{x} \notin \cup_{\iota \in \mathcal I_1 \cup \mathcal I_2} \mathcal D_\iota(\bF_p) \cup \mathcal E(\bF_p)$. Then it follows from Lemmas \ref{lemm:useful} and \ref{lemm:usefulII} that
\begin{align*}
\int_{\rho^{-1}(\tilde{x})} &= \left(1-\frac{1}{p}\right)^{-1}\int_{\mm_p^2}|y|_p^{s_\iota -d_\iota} \overline{\psi}_p(u/y^{e_\iota}) \,\mathrm dy_p \mathrm dz_p\\
&= \frac{1}{p-1} \int_{\mm_p}|y|_p^{s_\iota -d_\iota} \int_{\bZ_p^{\times}} \overline{\psi}_p(ub^{e_\iota}/y^{e_\iota})\, \mathrm db_p^{\times} \mathrm dy_p\\
&=
\begin{cases}
0 & \text{ if } e_\iota >1\\
-\frac{p^{-(s_\iota - d_\iota +2)}}{p-1} & \text{ if } e_\iota =1.
\end{cases}
\end{align*}
If $\tilde{x} \in \mathcal D_\iota(\bF_p) \cap \mathcal E(\bF_p)$ for some $\iota \in \mathcal I_3$, then we have
\begin{align*}
\int_{\rho^{-1}(\tilde{x})} &= \left(1-\frac{1}{p}\right)^{-1} \int_{\mm_p^2}|y|_p^{s_\iota -d_\iota} \overline{\psi}_p(z/y^{e_\iota}) \mathrm dy_p \mathrm dz_p\\
&= \left(1-\frac{1}{p}\right)^{-1} \int_{m_p}|y|_p^{s_\iota -d_\iota} \int_{\mm_p} \overline{\psi}_p(z/y^{e_\iota})\, \mathrm dz_p  \mathrm dy_p\\
&=
\begin{cases}
0 & \text{ if } e_\iota > 1\\
p^{-(s_\iota - d_\iota +2)} & \text{ if } e_\iota =1.
\end{cases}
\end{align*}
If $\tilde{x} \in \mathcal D_\iota(\bF_p) \cap \mathcal D_{\iota'}(\bF_p)$ for some $\iota, \iota' \in \mathcal I_3$, then it follows from Lemma \ref{lemm:usefulII} that
\begin{align*}
\int_{\rho^{-1}(\tilde{x})} &= \left(1-\frac{1}{p}\right)^{-1}\int_{\mm_p^2} |y|_p^{s_\iota -d_\iota} |z|_p^{s_{\iota'} -d_{\iota'}} \overline{\psi}_p \left( \frac{u}{y^{e_\iota} z^{e_{\iota'}}} \right)\, \mathrm dy_p \mathrm dz_p\\
&= \left(1-\frac{1}{p}\right)^{-1}\int_{\mm_p^2} |y|_p^{s_\iota -d_\iota} |z|_p^{s_{\iota'} -d_{\iota'}} \int_{\bZ_p^{\times}} \overline{\psi}_p \left( \frac{ub^{e_\iota}}{y^{e_\iota} z^{e_{\iota'}}} \right) \mathrm db_p^{\times} \, \mathrm dy_p \mathrm dz_p\\
&= 0.
\end{align*}
Thus we can conclude from these estimates and Lemma 9.4 in \cite{VecIII} that there exists an uniform bound $C'>0$ such that
$$
\left| \widehat{\sH}_p(\mathbf{s}, \lambda, t, \alpha) - 1 \right| < \frac{C'}{p^{1+\epsilon}}
$$
Our assertion follows from this.
\end{proof}

\begin{lemm}
\label{lemm:bounds}
Let $\Omega_\epsilon'$ be an open convex cone, defined by
$$
\begin{cases}
s_\iota - d_\iota -2 - \epsilon > 0  & \textnormal{ if } \iota \in \mathcal I_1\\
s_\iota - d_\iota +2 + 2\epsilon> 0 & \textnormal{ if } \iota \in \mathcal I_2\\
s_\iota -d_\iota +1> 0 & \textnormal{ if } \iota \in \mathcal I_3\\
\end{cases}
$$
where $\epsilon > 0$ is sufficiently small. Fix a compact subset of $\Omega_\epsilon'$ and $\epsilon \gg \delta > 0$. Then there exists a constant $C>0$ such that
$$
|\prod_p \widehat{\sH}_p(\mathbf{s}, \lambda, t, \alpha) \cdot \widehat{\sH}_\infty (\mathbf{s}, \alpha, t)| < \frac{C}{(1+t^2)|\beta|^{\frac{4}{3} - \frac{8}{3}\epsilon -\delta} |\gamma|^{1+\epsilon -\delta}},
$$
for $\Re(s)$ in that compact set, where $\alpha = \frac{\beta}{\gamma}$ with $\gcd (\beta, \gamma) =1$.
\end{lemm}
\begin{proof}
This lemma follows from Lemmas \ref{lemm:denominator}, \ref{lemm:numerator}, and \ref{lemm:real}, and from the proof of Lemma \ref{lemm:Euler_product}.
\end{proof}

\begin{theo}
\label{theo:main}
The zeta function $\sZ_1(\mathbf s, \textnormal{id})$ is holomorphic on the tube domain over an open neighborhood of the shifted effective cone $-K_X + \Lambda_{\textnormal{eff}}(X)$.
\end{theo}
\begin{proof}
Let $1 \gg \epsilon \gg \delta > 0$. Lemma \ref{lemm:bounds} implies that
$$
\sZ_1(\mathbf s, \text{id}) = \sum_{\lambda \in \mathbf M, \, \lambda(-1) = 1} \sum_{\alpha \in \bQ^{\times}} \frac{1}{2\pi} \int_{-\infty}^{+\infty} \prod_p \widehat{\sH}_p (\mathbf s, \lambda, t, \alpha) \cdot \widehat{\sH}_{\infty}(\mathbf s, t, \alpha)\, \mathrm dt,
$$
is absolutely and uniformly convergent on $\Omega_\epsilon'$, so $\sZ_1(\mathbf s, \text{id})$ is holomorphic on $\mathsf T_{\Omega_\epsilon'}$. Now note that the image of $\Omega_\epsilon'$ by $\Pic^G(X) \ra \Pic (X)$ contains an open neighborhood of $-K_X + \Lambda_{\textnormal{eff}}(X)$. Thus our theorem is concluded.
\end{proof}

\bibliographystyle{alpha}
\bibliography{axb}

\end{document}